\documentclass[12pt]{amsart}

\usepackage{amssymb}
\usepackage{hyperref} 
\usepackage{enumitem}

\usepackage{amsmath}
\usepackage{latexsym}
\usepackage{extarrows}
\usepackage{enumerate}
\usepackage{txfonts}
\usepackage{mathtools}
\usepackage{bm}
\usepackage{tikz-cd}
\usepackage{xcolor}
\usepackage{quiver}

\makeatletter
\@namedef{subjclassname@2020}{%
  \textup{2020} Mathematics Subject Classification}
\makeatother

\usepackage[T1]{fontenc}

\theoremstyle{plain}

\newtheorem{theorem}{Theorem}[section]
\newtheorem{proposition}[theorem]{Proposition}

\newtheorem{lemma}[theorem]{Lemma}

\theoremstyle{definition}

\newtheorem{definition}[theorem]{Definition}

\newtheorem{remark}[theorem]{Remark}

\newcommand\E{\mathbb{E}}
\newcommand\Z{\mathbb{Z}}
\newcommand\R{\mathbb{R}}
\newcommand\T{\mathbb{T}}
\newcommand\C{\mathbb{C}}

\newcommand\N{\mathbb{N}}
\newcommand\F{\mathbb{F}}

\newcommand\X{\mathcal{X}}
\newcommand\Y{\mathcal{Y}}

\newcommand\Hom{\operatorname{Hom}}
\newcommand\Aut{\operatorname{Aut}}
\newcommand\eps{\varepsilon}
\newcommand\Stone{\mathtt{Conc}}
\newcommand\CHProb{{\mathbf{CHPrb}}}
\newcommand\ProbAlg{{ \mathbf{PrbAlg}}}
\newcommand\OpProbAlg{{\mathbf{PrbAlg}}}
\newcommand\OpProbAlgG{{\mathbf{PrbAlg}_\Gamma}}
\newcommand\ConcProb{\mathbf{CncPrb}}
\newcommand\Cat{\mathcal{C}}

\begin{document}


\baselineskip=17pt


\title[]{The structure of arbitrary Conze--Lesigne systems}

\author[A. Jamneshan]{Asgar Jamneshan}
\address{Department of Mathematics\\ Ko\c{c} University \\ 
Rumelifeneri Yolu \\
34450, Sariyer, Istanbul, Turkey}
\email{ajamneshan@ku.edu.tr}

\author[O. Shalom]{Or Shalom}
\address{Einstein Institute of Mathematics \\ The Hebrew University of Jerusalem\\ 
Jerusalem \\
9190401, Israel}
\email{Or.Shalom@mail.huji.ac.il}

\author[T. Tao]{Terence Tao}
\address{Department of Mathematics\\ University of California \\ 
Los Angeles \\
CA 90095-1555, USA}
\email{tao@math.ucla.edu}

\date{\today}

\begin{abstract}  Let $\Gamma$ be a countable abelian group.  An (abstract) $\Gamma$-system $\mathrm{X}$ - that is, an (abstract) probability space equipped with an (abstract) probability-preserving action of $\Gamma$ - is said to be a \emph{Conze--Lesigne system} if it is equal to its second Host--Kra--Ziegler factor $\mathrm{Z}^2(\mathrm{X})$.  The main result of this paper is a structural description of such Conze--Lesigne systems for arbitrary countable abelian $\Gamma$, namely that they are the inverse limit of translational systems $G_n/\Lambda_n$ arising from locally compact nilpotent groups $G_n$ of nilpotency class $2$, quotiented by a lattice $\Lambda_n$.  Results of this type were previously known when $\Gamma$ was finitely generated, or the product of cyclic groups of prime order.  In a companion paper, two of us will apply this structure theorem to obtain an inverse theorem for the Gowers $U^3(G)$ norm for arbitrary finite abelian groups $G$.
\end{abstract}



\maketitle

\section{Introduction}

Furstenberg's ergodic-theoretic proof in \cite{furstenberg1977ergodic} of Szemer\'edi's theorem \cite{szemeredi1975sets} pioneered an influential synergy between ergodic theory and arithmetic combinatorics that continues to thrive in contemporary mathematics. Szemer\'edi's theorem asserts that a set $E\subset \Z$ such that $$\limsup_{N\to \infty} \frac{|E\cap\{-N,\ldots,N\}|}{2N+1}>0$$ contains arbitrarily long arithmetic progressions. By establishing a correspondence principle, Furstenberg showed that Szemer\'edi's theorem is equivalent to the multiple recurrence asymptotics  
\begin{equation}\label{multiple-recurrence}
\liminf_{N\to \infty} \frac{1}{2N+1} \sum_{n=-N}^N \mu(A\cap T^n(A)\cap T^{2n}(A) \cap \ldots \cap T^{(k-1)n}(A))>0
\end{equation}
for any measure-preserving transformation $T$ on a probability space $(X,\mu)$ and every subset $A$ in $X$ with $\mu(A)>0$. Furstenberg's novel approach laid the cornerstone for the development of ergodic Ramsey theory, giving rise to a series of intricate extensions of Szemer\'edi's theorem. These extensions include the multidimensional Szemer\'edi theorem \cite{furstenberg1978ergodic}, the density Hales--Jewitt theorem \cite{furstenberg1991density}, and the polynomial Szemer\'edi theorem \cite{bergelson1996polynomial}. It is noteworthy that alternative proofs for these extensions were discovered only much later. 

From an ergodic theoretic perspective, Furstenberg's multiple recurrence theorem stands as a significant extension of Poincar\'e's single recurrence theorem \cite{poincare}. While Poincar\'e's theorem can be proven succinctly and directly in modern treatments, Furstenberg's theorem is based on a structure theory, as initially developed by Furstenberg in \cite{furstenberg1977ergodic} and concurrently by Zimmer in \cite{zimmer1976ergodic,zimmer}. At the heart of this theory lies the Furstenberg--Zimmer structure theorem, which decomposes any measure-preserving system into a chain of more well-behaved subsystems.  

Nonconventional ergodic averages 
\begin{equation}\label{nonconventional-averages}
\frac{1}{2N+1} \sum_{n=-N}^N \prod_{i=1}^{k} f_i\circ T^{in}    
\end{equation}
serve as the functional counterparts to the multiple recurrence events in equation \eqref{multiple-recurrence}. This parallels how the conventional averages $\frac{1}{2N+1} \sum_{n=-N}^N f\circ T^n$ in von Neumann's mean ergodic theorem \cite{vonneumann} are the functional analog of single recurrence events in Poincar\'e's recurrence theorem. Von Neumann's mean ergodic theorem establishes the $L^2$-limit of these conventional averages and characterizes this limit as the projection onto the subspace of $T$-invariant functions. For a significant duration, an open question pertained to the existence of an analogue of von Neumann's mean ergodic theorem for nonconventional ergodic averages.

Partial progress in this direction was obtained by Conze and Lesigne in a series of papers \cite{cl1,cl2,cl3}. They proved the convergence of the nonconventional ergodic averages \eqref{nonconventional-averages} in the particular case $k=3$ (under the additional assumption that the underlying system is totally ergodic). Moreover, they characterized the $L^2$-limit by identifying a subsystem of the underlying $\Z$-systems which is characteristic for these averages, that is, the nonconventional average of the projection of the functions $f_i$ on the $L^2$-subspace of  this subsystem has the same limit as the nonconventional average of the original functions $f_i$.

Much later, a complete answer was provided independently by breakthrough works of Host and Kra \cite{host2005nonconventional} and Ziegler \cite{ziegler2007universal}.  
They confirmed the $L^2$-convergence of the nonconventional ergodic averages \eqref{nonconventional-averages} by a significant refinement of the Furstenberg--Zimmer structure theory that identified a finer hierarchy of subsystems  of the original system which effectively control or are characteristic for these averages for all $k$. In this refined Host--Kra--Ziegler tower, each $k$-th subsystem is referred to as the "factor of order $k$", while the order $2$ factor corresponds to the subsystem identified by Conze and Lesigne. Their deep insight lies in an algebraic and geometric classification of these order $k$ factors, identifying them as inverse limits of systems formed by translations on nilmanifolds of the form $G_n/\Lambda_n$, where $G_n$ is a nilpotent Lie group of nilpotency class $k$, and $\Lambda_n\leq G_n$ denotes a lattice - a discrete and cocompact subgroup of $G_n$. 

In a parallel development, Gowers introduced a new influential Fourier-analytic proof of Szemer\'edi's theorem in \cite{gowers1,gowers2} that represents a substantial generalization of Roth's approach \cite{roth1952} in the case of $3$-term progressions. It marked the inception of the field of higher-order Fourier analysis, with the Gowers $U^k$ norms on finite abelian groups $G$ taking center stage. The $U^k$ norm essentially quantifies the normalized average of $2^k$-fold autocorrelations of functions $f\colon G\to \C$ over arithmetic cubes of the form $(x+\omega\cdot h)_{\omega\in \{0,1\}^k}$ for $h\in G^{\{0,1\}^k}$. A crucial aspect of the Gowers norms within additive combinatorics and their applications to analytic number theory lies in their control over certain "bounded-complexity" multilinear forms on finite abelian groups such as forms associated with arithmetic progressions. To harness the control over these multilinear forms, one must solve the inverse problem for the Gowers uniformity norms, which, roughly speaking, asks for an algebraic classification of $1$-bounded functions $f\colon G\to\C$ with a large $U^k$ norm for all $k\geq 3$  (where the cases of $k=1$ and $k=2$ can be readily derived from basic Fourier analysis). A foundational achievement in higher-order Fourier analysis credited to Green, the third author, and Ziegler in their work \cite{gtz} is the resolution of the inverse problem concerning the Gowers $U^k$ norms for cyclic groups. Their inverse theorem asserts that a function with positive Gowers $U^k$ norm exhibits correlation with a function derived from a nilmanifold (having nilpotency class $k$). The Green--Tao--Ziegler inverse theorem has important applications in analytic number theory such as establishing the correct asymptotics for primes in arithmetic progressions  \cite{gt-linear}. Despite significant progress in important special cases, an inverse theorem for the Gowers $U^k$ norms on arbitrary finite abelian groups is currently open for all $k\geq 4$. The $k=3$ case was resolved by two of us in the companion paper \cite{jt21-1}.   

The heuristic analogy between the Host--Kra--Ziegler structure theorem and the Green--Tao--Ziegler inverse theorem is remarkably compelling. There is substantial support for this analogy. Host and Kra characterized the order $k$ factors of $\mathbb{Z}$-systems by the Host--Kra--Gowers seminorms of order $k+1$, which serve as an infinitary counterpart to the Gowers uniformity norms. 
In this context, the Host--Kra--Ziegler structure theorem can be viewed as a resolution of the inverse problem associated with these Host--Kra--Gowers seminorms for $\Z$-actions. 

The present paper and its companion paper \cite{jt21-1} contribute to these developments, specifically enhancing our comprehension of the above heuristic analogy. In this paper, we establish a Host--Kra-type structure theorem for arbitrary abelian systems of order two, also called \emph{Conze--Lesigne systems}.  
Questions of recurrence and convergence of nonconventional ergodic averages are not addressed in this paper. However, the results of this paper will be applied in the companion paper \cite{jt21-1} to give a qualitative proof of the inverse theorem for the Gowers uniformity norms $U^3(G)$ for arbitrary finite abelian groups $G$ via a correspondence principle. 

\subsection{A note on probability space conventions}\label{note-sec}

For technical reasons we will have to work in this paper with three slightly different categories of probability spaces, as well as their corresponding categories of measure-preserving systems associated to a group $\Gamma$ (which will be countable\footnote{In this paper we use ``countable'' as an abbreviation for ``at most countable''.} in most, though not all, of the contexts we will consider):

\begin{itemize}
    \item[(i)] The category of \emph{concrete probability spaces} $(X,\X,\mu)$, in which one can meaningfully talk about individual points $x$ in the space $X$, and maps between these spaces are defined in a pointwise fashion. One can then form the category of \emph{concrete $\Gamma$-systems} $(X,\X,\mu,T)$ of concrete probability spaces equipped with a pointwise defined measure-preserving action $T \colon \gamma \mapsto T^\gamma$ of $\Gamma$.   Among other things, concrete probability spaces are a convenient category in which to study group extensions by measurable (but not necessarily continuous) cocycles.
    \item[(ii)] The category of \emph{probability algebras}\footnote{In the language of \cite[Chapter 2]{glasner2015ergodic}, probability algebras are referred to as \emph{measure algebras}, and (separable) abstract $\Gamma$-systems are referred to as \emph{measure-preserving dynamical systems}. In the language of \cite{jt-foundational}, probability algebras are essentially abstract probability spaces with the additional property that all non-empty abstract subsets have positive measure.} $(\X,\mu)$, in which one has ``quotiented out all the null sets''; as a consequence, one can no longer meaningfully refer to individual points, and maps between spaces are typically only defined up to almost everywhere equivalence.  One can then form the category of \emph{abstract $\Gamma$-systems} $(\X,\mu, T)$ of probability algebras equipped with an abstract measure-preserving action $T$ of $\Gamma$.  The category of abstract $\Gamma$-systems is the most natural category in which to discuss factors of a system, such as the Host--Kra--Ziegler or Conze--Lesigne factors, as well as to discuss the isomorphic nature of two systems.
    \item[(iii)] The category of \emph{compact probability spaces} $(X,{\mathcal F}, \X, \mu)$, in which the probability spaces are now compact Hausdorff (with the measure $\mu$ being a (Baire-)Radon measure), and the maps between spaces are now additionally required to be continuous.   This then forms the category of \emph{compact $\Gamma$-systems} $(X,{\mathcal F}, \X, \mu, T)$ of compact probability spaces equipped with a \emph{continuous} measure-preserving action $T$ of $\Gamma$. The category of compact $\Gamma$-systems is the most natural category to discuss transitivity properties of a group action, or to compute the stabilizer of such an action at a point.
\end{itemize}

For the convenience of the reader we review the definition of these categories, as well as the relationships between them that we will need, in Appendix \ref{conc-app}. Very roughly speaking, as long as one is in the ``countable'' setting in which the acting group $\Gamma$ is countable, the concrete probability spaces are Lebesgue spaces, the probability algebras are separable, and the compact Hausdorff probability spaces are metrizable, then these three categories are ``morally interchangeable'', largely thanks to the ability to construct \emph{topological models} of abstract $\Gamma$-systems (and \emph{continuous representations} of factor maps between such systems); however more care needs to be taken in ``uncountable'' settings when one or more of the above assumptions is not in force, and even in the countable setting there are some subtleties, particularly with regard\footnote{See in particular the erratum \cite{hk-errata} to \cite[Chapter 19]{hk-book} for further discussion of this particular subtlety.}
to ``near-actions'' on concrete probability spaces that are only defined up to almost everywhere equivalence (see Appendix \ref{conc-app} for the definition of near-actions), and are thus only genuine actions in an abstract sense.  Most of our arguments will take place within a countable setting (and are already new in this case), but through appropriate use of inverse limits our main result will also be applicable for inseparable\footnote{Furthermore, even if one is only interested in applying our results for separable systems, there is one step in the argument in which a potentially inseparable system can arise, namely when one uses Gelfand duality to construct a (possibly inseparable) topological model of a (separable) system, which we call a \emph{Koopman model}; see Appendix \ref{koopman-sec}.  While it is possible with significant further effort to demonstrate that this model is in fact separable in the specific context being considered (cf., \cite{hk-errata}), it shortens the arguments to just proceed without verifying separability, as this property turns out to not be needed in the proofs.}  systems.  For a first reading we recommend that the reader ignore the fine technical distinctions between these categories, or between the countable and uncountable cases.

Traditionally, the literature has been focused on concrete Lebesgue $\Gamma$-systems.  However, it will be convenient to phrase our main results in the setting of abstract (and not necessarily separable) $\Gamma$-systems, although thanks to the aforementioned modeling results one can often reformulate these results in the other categories mentioned, particularly in the separable case.  In particular, the factor relation $\mathrm{Y} \leq \mathrm{X}$ between two $\Gamma$-systems $\mathrm{X}, \mathrm{Y}$ (as defined in Appendix \ref{conc-app}) will be understood to be in the abstract sense unless otherwise specified, even when the systems $\mathrm{X}, \mathrm{Y}$ can be viewed as concrete or compact $\Gamma$-systems; similarly for the notion of an inverse limit of factors.  While the factor maps $\pi \colon \mathrm{X} \to \mathrm{Y}$ in this paper are initially only defined abstractly, in practice they can often be upgraded to concrete measurable maps by using tools such as those in Proposition \ref{reverse}.

\begin{remark}
If we restrict ourselves to separable probability algebras, it may be possible to replace the category of probability algebra dynamical systems with the category of Lebesgue probability spaces equipped with near-actions. However, in order to do this, it is necessary to implement estimates similar to those established in \cite{hk-errata} to ensure that near-actions of Polish groups can be accurately described by a continuous action on a \emph{separable} model.  
\end{remark}

We isolate some key special examples of compact $\Gamma$-systems (which can then be viewed as concrete or abstract $\Gamma$-systems by forgetting some of the structure):

\begin{definition}[Translational and rotational $\Gamma$-systems]  Let $\Gamma$ be a group.
A \emph{translational $\Gamma$-system} is a compact $\Gamma$-system of the form $G/\Lambda = (G/\Lambda, \mu, T)$, where $G = (G,\cdot)$ is a locally compact\footnote{In this paper we use ``locally compact group'' as shorthand for ''locally compact Hausdorff second countable group''.  Similarly for ``compact group'' or ``compact abelian group''.} unimodular\footnote{The unimodularity hypothesis is required in order to have a well-defined Haar measure on the quotient space $G/\Lambda$.  In our applications, the locally compact group $G$ will be nilpotent and thus automatically unimodular.} group, $\Lambda$ is a closed cocompact subgroup of $G$, $\mu$ is the Haar probability measure on the compact quotient space $G/\Lambda$, and the action $T$ is given by $T^\gamma x = \phi(\gamma) x$ for all $\gamma \in \Gamma$, $x \in G/\Lambda$, for some homomorphism $\phi \colon \Gamma \to G$.  If $G$ is a compact abelian group (which we now write additively as $Z=(Z,+)$) and $\Lambda$ is trivial, we refer to the translational $\Gamma$-system $Z = (Z, \mu, T)$ as a \emph{rotational $\Gamma$-system}.
\end{definition}

Amongst the translational $\Gamma$-systems $G/\Lambda$, we single out for special mention the \emph{$\Gamma$-nilsystems} of order at most $k$ for a given $k \geq 1$, in which $G$ is a nilpotent Lie group of nilpotency class at most $k$, and $\Lambda$ a lattice (i.e., a discrete cocompact subgroup) of $G$.  For instance, rotational $\Gamma$-systems are inverse limits of $\Gamma$-nilsystems of order at most $1$.

\subsection{Host--Kra--Ziegler factors and Conze--Lesigne systems} 

Let $\Gamma = (\Gamma,+)$ be a countable discrete abelian group, and let $\mathrm{X} = (\X,\mu, T)$ be an (abstract) $\Gamma$-system.  Among the factors of $\mathrm{X}$ we can form the \emph{invariant factor} $\mathrm{Z}^0(\mathrm{X})$, defined by replacing the $\sigma$-complete Boolean algebra $\X$ by its invariant subalgebra
$$ \X^T  \coloneqq \bigcap_{\gamma \in \Gamma} \{ E \in \X: E = (T^\gamma)^* E \}$$
and restricting $\mu$ and $T$ accordingly.  As usual, we call the $\Gamma$-system $\mathrm{X}$ \emph{ergodic} if this invariant factor is trivial.  
Similar notions can now be defined for concrete or compact $\Gamma$-systems by forgetting some of the structure.  For most of this paper we will focus on ergodic systems; in principle one can use tools such as the ergodic decomposition (or conditional analysis, see \cite{jamneshan2018measure}) to adapt the results in this paper to the non-ergodic setting, but we will not attempt to do so here.

The invariant factor $\mathrm{Z}^0(\mathrm{X})$ of an (abstract) $\Gamma$-system $\mathrm{X}$ is the zeroth in the sequence of \emph{Host--Kra--Ziegler factors}
$$\mathrm{Z}^0(\mathrm{X}) \leq \mathrm{Z}^1(\mathrm{X}) \leq \mathrm{Z}^2(\mathrm{X}) \leq \dots \leq \mathrm{X}$$
of $\mathrm{X}$; we briefly review the precise definition of these factors in Section \ref{notation-sec}.  We will not directly use this definition as we will rely on existing results about these factors in the literature, but we will remark that $\mathrm{Z}^k(\mathrm{X})$ is the universal characteristic factor for the \emph{Host--Kra--Gowers seminorm} $\| \cdot \|_{U^{k+1}(\mathrm{X})}$ on $\mathrm{X}$; see e.g., \cite{host2005nonconventional}, \cite[Appendix A]{btz}, \cite{hk-book}.  These norms are traditionally defined for concrete Lebesgue $\Gamma$-systems, but their definitions can be easily adapted to the abstract setting, or alternatively one can replace an abstract $\Gamma$-system by a suitable concrete (or topological) model and apply the standard constructions to that model; see Section \ref{notation-sec}.

The first Host--Kra--Ziegler factor $\mathrm{Z}^1(\mathrm{X})$ is known as the \emph{Kronecker factor} and was studied by von Neumann and Halmos, see, e.g., \cite{halmos2013lectures} for a reference.  The second Host--Kra--Ziegler factor $\mathrm{Z}^2(\mathrm{X})$ is known as the \emph{Conze--Lesigne factor} and was studied (in the $\Gamma=\Z$ case at least) by Conze and Lesigne \cite{cl2}, \cite{cl3} (see also \cite{rudolph}, \cite{meiri}, \cite{fw}, \cite{hk1}, \cite{hk2}).  For general $k$, the Host--Kra--Ziegler factors were introduced in the $\Gamma=\Z$ case by Host and Kra \cite{host2005nonconventional}; in the subsequent work of Ziegler \cite{ziegler2007universal} the universal characteristic factors for multiple recurrence were constructed, which were later shown by Leibman (see  \cite[Appendix A]{leibman-hkz} to be equivalent to the factors of Host and Kra.  As is well known, the constructions of the Host--Kra--Ziegler factors, the Host--Kra--Gowers norms, and the Host--Kra parallelepiped systems extend without difficulty to arbitrary countable abelian groups $\Gamma$; see for instance \cite[Appendix A]{btz}, where the factor $\mathrm{Z}^k(\mathrm{X})$ was denoted instead as $\mathrm{Z}_{<k+1}(\mathrm{X})$.

Let $k \geq 0$ be a natural number.  An ergodic (abstract) $\Gamma$-system $\mathrm{X}$ is said to be of \emph{order (at most) $k$} if $\mathrm{X} = \mathrm{Z}^k(\mathrm{X})$.  Thus for instance an ergodic $\Gamma$-system is of order $0$ if and only if it is (abstractly) trivial.  We recall some simple facts about such systems:

\begin{lemma}[Basic facts about systems of order $k$]\label{basic-facts}  Let $\Gamma$ be a countable abelian group.
\begin{itemize}
\item[(i)]  $\mathrm{Z}^k(\mathrm{X})$ is of order $k$ for any ergodic (abstract) $\Gamma$-system $\mathrm{X}$.
\item[(ii)] Any factor of an ergodic $\Gamma$-system of order $k$ will also be an ergodic $\Gamma$-system of order (at most) $k$.  
\item[(iii)]  The inverse limit of ergodic $\Gamma$-systems of order $k$ will also be an ergodic $\Gamma$-system of order $k$.
\end{itemize}
\end{lemma}
 
\begin{proof} For (i), see \cite[Corollary 4.4]{host2005nonconventional}, \cite[(A.9)]{btz}, or \cite[Chapter 9, Theorem 15(ii)]{hk-book}.  For (ii), see \cite[Proposition 4.6]{host2005nonconventional}, \cite[Lemma A.34]{btz}, or \cite[Chapter 9, Proposition 17(ii)]{hk-book}.  For (iii), see \cite[Proposition 4.6]{host2005nonconventional}, \cite[Lemma A.34]{btz}, or \cite[Chapter 9, Theorem 20]{hk-book}.
\end{proof}

Ergodic $\Gamma$-systems of order $1$ will be referred to as \emph{Kronecker systems}, while ergodic $\Gamma$-systems of order $2$ will be referred to as \emph{Conze--Lesigne systems}; thus for instance a Conze--Lesigne system is its own Conze--Lesigne factor.  The classification of systems of arbitrary order is of importance in the theory of multiple recurrence; for instance, as seen in \cite{host2005nonconventional},  \cite{ziegler2007universal}, classification of ergodic separable $\Z$-systems $(X,\X,\mu,T)$ of order $k$ was used to give the first proofs of the norm convergence of multiple ergodic averages $\frac{1}{N} \sum_{n=1}^N T^n f_1 \dots T^{(k+1) n} f_{k+1}$ for $f_1,\dots,f_{k+1} \in L^\infty(X,\X,\mu)$ for general $k$.

The classification of Kronecker systems is well-known (going back to the work of von Neumann and Halmos \cite{halmos2013lectures}):

\begin{theorem}[Classification of Kronecker systems]\label{kron}  Let $\Gamma$ be a countable abelian group and let $\mathrm{X}$ be an ergodic separable $\Gamma$-system.  Then the following are equivalent:
\begin{itemize}
\item[(i)] $\mathrm{X}$ is a Kronecker $\Gamma$-system.
\item[(ii)]  $\mathrm{X}$ is (abstractly) isomorphic to a rotational $\Gamma$-system $Z$ for some compact abelian metrizable group $Z$.
\item[(iii)] $\mathrm{X}$ is the inverse limit of rotational $\Gamma$-systems $Z_n$ for some compact abelian \emph{Lie} groups $Z_n$.
\end{itemize}
\end{theorem}

\begin{proof}  For $\Z$-systems, the equivalence of (i) and (ii) follows from \cite[Chapter 2, Theorem 12]{hk-book} and \cite[Chapter 9, Proposition 8]{hk-book} (see also the discussion before \cite[Lemma 4.2]{host2005nonconventional}); the arguments extend without difficulty to arbitrary countable abelian $\Gamma$ (and one can replace the abstract $\Gamma$-system by a concrete model if desired, or argue directly in the abstract setting).  The deduction of (i) from (iii) then follows from Lemma \ref{basic-facts}(iii), while the deduction of (iii) from (ii) follows from the Peter--Weyl theorem (or Pontryagin duality); see e.g., \cite[Exercise 1.4.26]{tao-hilbert}.
\end{proof}

We remark that it is not difficult to remove the separability and countability hypotheses from Theorem \ref{kron}, so long as one similarly removes the metrizability hypothesis from conclusion (ii).  As a consequence of this theorem (and Lemma \ref{basic-facts}), the Kronecker factor $\mathrm{Z}^1(\mathrm{X})$ of an ergodic concrete $\Gamma$-system $\mathrm{X}$ can be equivalently described as the maximal rotational factor of $\mathrm{X}$ (cf. \cite[Proposition 13(iv)]{hk-book}). 

Now we turn to the higher order Host--Kra--Ziegler factors.  In the case of $\Z$-systems, we have the following fundamental result of Host and Kra:

\begin{theorem}[Classification of Host--Kra--Ziegler $\Z$-systems]\label{hk-thm}
Let $k \geq 1$ be a natural number, and let $\mathrm{X}$ be an ergodic separable $\Z$-system.  Then the following are equivalent:
\begin{itemize}
\item[(i)] $\mathrm{X}$ is a $\Z$-system of order (at most) $k$.
\item[(ii)] $\mathrm{X}$ is the inverse limit of $\Gamma$-nilsystems $G_n/\Lambda_n$ of order at most $k$ (as defined at the end of Section \ref{note-sec}).
\end{itemize}
\end{theorem}

The implication of (i) from (ii) can be found in \cite[Chapter 12, Corollary 19]{hk-book}; the implication of (ii) from (i) is more difficult and was proven in \cite[Theorem 10.1]{host2005nonconventional} (see also \cite{ziegler2007universal} for closely related results, and \cite{hk-book} for a more detailed exposition). For a treatment of the Conze--Lesigne case $k=2$, see \cite{cl1,cl2,cl3}, \cite[\S 3]{meiri}, \cite{rudolph}, \cite{fw}, \cite[\S 8]{host2005nonconventional}, \cite[\S 9]{ziegler-thesis}.  An alternate proof of this theorem using compact nilspaces was also given in \cite{gutman-lian}.  

Theorem \ref{hk-thm} was extended to $\Z^d$-systems for any finite $d$ by Griesmer \cite[Theorem 4.1.2]{griesmer}, following similar arguments to those in \cite{host2005nonconventional}; a further extension to $\Gamma$-systems for any finitely generated nilpotent group $\Gamma$ (extending the preceding definitions to the nilpotent case in a natural fashion) was obtained using the machinery of nilspaces in \cite[Theorem 5.12]{cs-cubic} (with the abelian case previously established by this method in \cite{gutman-lian}).  However, the situation changes somewhat once one considers groups $\Gamma$ that are not finitely generated; in particular, the arguments in \cite{host2005nonconventional}, \cite{griesmer} rely crucially on finite generation to establish some connectedness properties of certain structure groups arising in the analysis that do not hold in general in the infinitely generated case.  A model infinitely generated case is that of the countably generated vector space $\F_p^\omega = \bigcup_{n=1}^\infty \F_p^n$ over a finite field $\F_p$ of prime order $p$. In this case we have a fairly satisfactory classification, particularly in the case of high characteristic:

\begin{theorem}[Classification of Host--Kra--Ziegler $\F_p^\omega$-systems]\label{btz-thm}
Let $k \geq 1$ be a natural number, let $p$ be a prime, and let $\mathrm{X}$ be an ergodic separable $\F_p^\omega$-system.  In the high characteristic case $p \geq k-1$, then the following are equivalent:
\begin{itemize}
\item[(i)] $\mathrm{X}$ is a $\F_p^\omega$-system of order (at most) $k$.
\item[(ii)] $\mathrm{X}$ is generated by phase polynomials\footnote{We will not need the concept of a phase polynomial to state or prove our main results, but see for instance  \cite[Definition 1.13]{btz} for a precise definition.} of degree at most $k$.  
\end{itemize}
In the low characteristic case $p<k-1$, (ii) still implies (i), but it is currently open whether (i) implies (ii) in these cases.  The weaker implication is known that ergodic separable $\F_p^\omega$-systems of order $k$ are generated by phase polynomials of degree at most $C(k)$ for some $C(k)$ depending only on $k$.
\end{theorem}

\begin{proof}  The implication of (ii) from (i) (in both high and low characteristic) is \cite[Lemma A.35]{btz}; the converse implication was established for $p \geq k+1$ in \cite[Theorem 1.18]{btz} and recently extended to $p=k, k-1$ in \cite[Theorem 1.12]{CGSS}.  The final claim of the theorem is \cite[Theorem 1.19]{btz}.
\end{proof}

The form of Theorem \ref{btz-thm} does not closely resemble that of Theorem \ref{hk-thm}.  In more recent work of the second author, results closer in appearance to Theorem \ref{hk-thm} were established for various classes of groups $\Gamma$:

\begin{theorem}[Partial classifications of Host--Kra--Ziegler systems]\label{partial}  Let $k \geq 1$, let $\Gamma$ be a countable abelian group.  Let $\mathrm{X}$ be an ergodic separable $\Gamma$-system of order $k$.
\begin{itemize}
    \item[(i)]  \cite[Theorem 1.31]{shalom3}  If $k=2$, and $\Gamma = \bigoplus_{p \in {\mathcal P}} \Z/p\Z$ for some countable multiset ${\mathcal P}$ of primes, then $\mathrm{X}$ is the inverse limit of translational systems $G_n/\Lambda_n$, where each $G_n$ is a locally compact Polish $2$-step nilpotent group.
    \item[(ii)] \cite[Theorem 2.3]{shalom1}  If $\Gamma = \F_p^\omega$ with $k \leq p-1$, then $\mathrm{X}$ is equivalent to a translational $\Gamma$-system $G/\Lambda$ with $G$ and $\Lambda$ totally disconnected and nilpotent of class at most $k$.  
    \item[(iii)] \cite[Theorem 2.10]{shalom1} If $\Gamma = \bigoplus_{p \in {\mathcal P}} \Z/p\Z$ for some countable multiset ${\mathcal P}$ of primes, then there exists a natural number $m = m(k)$ depending only on $k$, and an $m$-extension\footnote{See \cite[Definition 2.4]{shalom1} for a definition of this term, which we will not need in the rest of this paper.} $\mathrm{Y}$ of $\mathrm{X}$, which is an ergodic separable $\Gamma'$-system for some countable abelian group $\Gamma'$ which is the inverse limit of translational $\Gamma'$-systems $G_n/\Lambda_n$, where each $G_n$ is a finite dimensional\footnote{See \cite[Definition 2.6]{shalom1} for the definition of finite dimensionality for locally compact groups; we will not need this notion in the rest of this paper.} locally compact group of nilpotency class at most $k$, and $\Lambda_n$ is totally disconnected.
\end{itemize}
\end{theorem}

We also mention some further variants of the above results:
\begin{itemize}
    \item In \cite[Theorem 1.21]{shalom2}, the second author proved that when $\Gamma$ is countable abelian and $\mathrm{X}$ is a Conze--Lesigne $\Gamma$-system, there exists a nilpotent locally compact Polish group $G$, a compact totally disconnected group $K$, and a closed totally disconnected subgroup $\Lambda$ of $G$ such that $\mathrm{X}$ is (abstractly) isomorphic to the double coset system $K\backslash G/\Lambda$ acting by a translation action $T^\gamma x = \phi(\gamma) x$ for some homomorphism $\phi \colon \Gamma \to G$ that normalizes $K$. 
    \item In \cite[Theorem 1.18]{shalom2}, the second author showed that for any countable abelian group $\Gamma$ and ergodic separable $\Gamma$-system $\mathrm{X}$, there is an extension $\mathrm{Y}$ of $\mathrm{X}$ whose Conze--Lesigne factor is a translational $\Gamma'$-system $G/\Lambda$ for some extension $\Gamma'$ of $\Gamma$ and some locally compact Polish group $G$ that is nilpotent of nilpotency class at most two, where the notion of extension was defined in \cite{shalom2}.  
    \item In \cite[Theorem 1.9]{CGSS}, it was shown in both high and low characteristic that an ergodic separable $\F_p^\omega$-system of order $k$ is a $p$-homogeneous $k$-step nilspace system (see \cite{CGSS} for the definitions of these terms, which we will not need in the rest of this paper).
\end{itemize}

\subsection{Main result}

We now come to the main new result of this paper, which is to establish a complete description of Conze--Lesigne factors for arbitrary countable abelian groups $\Gamma$:

\begin{theorem}[Classification of Conze--Lesigne $\Gamma$-systems]\label{main-thm}
Let $\Gamma$ be a countable abelian group, and let $\mathrm{X}$ be an ergodic (abstract) $\Gamma$-system.  Then the following are equivalent:
\begin{itemize}
\item[(i)] $\mathrm{X}$ is a Conze--Lesigne $\Gamma$-system (i.e., a $\Gamma$-system of order at most $2$).
\item[(ii)] $\mathrm{X}$ is the inverse limit of translational  $\Gamma$-systems $G_n/\Lambda_n$, where each $G_n$ is a locally compact nilpotent Polish group of nilpotency class two, and $\Lambda_n$ is a lattice (i.e., a discrete cocompact subgroup) in $G_n$.  Furthermore, $G_n$ contains a closed central subgroup $G_{n,2}$ containing the commutator group $[G_n,G_n]$, with $\Lambda_n \cap G_{n,2}$ a lattice in $G_{n,2}$.
\end{itemize}
In (ii) we can also require that $G_{n,2}$ is a compact abelian Lie group, $\Lambda_n \cap G_{n,2}$ is trivial, and $\Lambda_n$ is abelian.
\end{theorem}

We remark that as $\Lambda_n$ is a discrete subgroup of the Polish group $G_n$, it is automatically countable.

This result can be compared with the previously mentioned result in \cite[Theorem 1.18]{shalom2}.  On the one hand, Theorem \ref{main-thm} does not require the passage to some extension $\mathrm{Y}$ of $\mathrm{X}$; on the other hand, the conclusion is weaker as the system is described as an inverse limit of translational systems, rather than as a translational system.  In view of Theorem \ref{hk-thm}, one could ask whether one could strengthen Theorem \ref{main-thm} further by requiring in (ii) that the $G_n$ are nilpotent Lie groups, and $\Lambda_n$ lattices, so that $\mathrm{X}$ would be the inverse limit of nilsystems.  Unfortunately when $\Gamma$ is not finitely generated, there are counterexamples that show that this stronger version of Theorem \ref{main-thm} fails; see the example presented after \cite[Conjecture 2.14]{shalom1} (in the discussion of \cite[Theorem 4.3]{shalom1}).  In Theorem \ref{main-thm} it is not required that the system $\mathrm{X}$  be separable, but it turns out it is quite easy to reduce to this case, and indeed this will be one of the first steps in the proof.  The group $G_{n,2}$ in Theorem \ref{main-thm} can in fact be taken to be the commutator group $[G_n,G_n]$ if desired; see Remark \ref{kron-remark}.  However, from the theory of filtered nilpotent groups (see e.g., \cite[Appendix B]{gtz}) it seems more natural to allow $G_{n,2}$ to be slightly larger than the commutator group (see Section \ref{first-ex} for one example of this).

In Section \ref{examples-sec} we provide some examples of Conze--Lesigne systems associated to groups in even, odd, and zero characteristic that illustrate the conclusion of Theorem \ref{main-thm} despite not being obviously associated to any nilpotent structures. 

\begin{remark}  It is tempting to speculate as to whether Theorem \ref{main-thm} can be extended to systems of order $k$ for $k>2$, by some induction on $k$.  Here one runs into a significant technical obstacle even when $k=3$; whereas in the $k=2$ case, the system can be expressed (using Theorem \ref{cl-extend} below) as a group extension of a translational system (indeed a rotational system, in this case), the analogous arguments in the $k=3$ case (when combined with Theorem \ref{main-thm}) only allow one to describe systems of order $3$ as group extensions of \emph{inverse limits} of translational systems.  When the group $\Gamma$ is finitely generated, one can use the connectedness of the structure groups to avoid this issue (cf. \cite[Lemma 10.4]{host2005nonconventional}), but in the infinitely generated case it is not clear to us whether such group extensions of inverse limits of translational systems can necessarily be expressed as inverse limits of translational systems, even if one possesses suitable higher-order analogues of the Conze--Lesigne equation.  We hope to investigate these issues further in subsequent work.
\end{remark}

\begin{remark}  Theorem \ref{main-thm} does not immediately imply previous structural results about Conze--Lesigne systems for specific groups $\Gamma$, such as those stated in Theorems \ref{hk-thm}, \ref{btz-thm}, \ref{partial}, because these theorems can take advantage of special features of the groups $\Gamma$ they consider to obtain stronger conclusions than those in Theorem \ref{main-thm}(ii), see Section \ref{sec-speculation} below for a related discussion.   However, one could certainly use Theorem \ref{main-thm} as a ``black box'' to \emph{shorten} the proofs of these other theorems, by reducing matters to the study of nilpotent translational systems $G/\Lambda$ of nilpotency class two, for which many of the intermediate statements used in the course of those proofs are easy to establish.  We leave the details of such shortenings to the interested reader. 
\end{remark}

\subsection{An application to the Gowers uniformity norms}

In \cite{taozieglerhigh}, the description of systems of order $k$ for $\F_p^\omega$-systems from Theorem \ref{btz-thm} was combined with a correspondence principle to establish an inverse theorem for the Gowers uniformity norm $U^{k+1}(\F_p^n)$ for finite-dimensional vector spaces $\F_p^n$.  In a similar vein, Theorem \ref{main-thm} can be combined with a correspondence principle to establish an inverse theorem for the Gowers norm $U^3(G)$ associated to an arbitrary finite abelian group $G$.  More precisely, in the companion paper \cite{jt21-1} to this paper, we show

\begin{theorem}[Inverse theorem for $U^3(G)$]\label{uk-inverse-nil}  Let $G$ be a finite additive group, let $\eta > 0$, and let $f \colon G \to \C$ be a $1$-bounded function with $\|f\|_{U^{3}(G)} \geq \eta$.  Then there exists a degree $2$ filtered nilmanifold $H/\Gamma$, drawn from some finite collection ${\mathcal N}_{\eta}$ of such nilmanifolds that depends only on $\eta$ but not on $G$ (and each such nilmanifold in ${\mathcal N}_{\eta}$ endowed arbitrarily with a smooth Riemannian metric), a Lipschitz function $F \colon H/\Gamma \to \C$ of Lipschitz norm $O_{\eta}(1)$, and a polynomial map $g \colon G \to H/\Gamma$ such that
\begin{equation}\label{fF}
 |\E_{x \in G} f(x) \overline{F(g(x))}| \gg_{\eta} 1.
 \end{equation}
\end{theorem}

We refer the reader to \cite{jt21-1} for definitions of all the terms in this theorem, as well as for details of how it is derived from Theorem \ref{main-thm}.

\subsection{Overview of proof}

Our proof of Theorem \ref{main-thm} is based primarily on the methods of Host and Kra \cite{host2005nonconventional}, \cite{hk-book}, while also incorporating some tools from \cite{btz}.  It is not difficult to reduce to the case of separable $\Gamma$-systems.  The next step, which is quite standard, is to express Conze--Lesigne systems as abelian extensions of the Kronecker factor. 

\begin{theorem}[Conze--Lesigne systems are abelian extensions of Kronecker factor]\label{cl-extend}  Let $\Gamma$ be a countable abelian group, and let $\mathrm{X}$ be an ergodic separable $\Gamma$-system.  Then the following are equivalent:
\begin{itemize}
    \item[(i)]  $\mathrm{X}$ is a Conze--Lesigne system.
    \item[(ii)] $\mathrm{X}$ is (abstractly) isomorphic to a group extension $Z \rtimes_\rho K$, where $Z$ is a rotational $\Gamma$-system, $K$ is a compact abelian group, and $\rho$ is a $(Z, K)$-cocycle of type $2$.
\end{itemize}
Furthermore, in part (ii), we can take $Z$ to be equivalent to the Kronecker factor.
\end{theorem}

The notions of cocycle, extension, and cocycle type appearing in the above theorem will be reviewed in Section \ref{cocyc} below.

\begin{proof} For $\Z$-systems, the implication of (ii) from (i) was established in \cite[Proposition 6.34]{host2005nonconventional} or \cite[Chapter 18, Theorem 6]{hk-book}; as observed previously by several authors \cite[Proposition 3.4]{btz}, \cite[Proposition 1.16]{shalom3}, \cite[Proposition A.18]{shalom1} (see also \cite[Proposition 3.6]{tz-concat} for an alternate proof of the abelian nature of $K$), the arguments extend without difficulty to arbitrary countable abelian groups $\Gamma$. The implication of (i) from (ii) follows from \cite[Corollary 7.7]{host2005nonconventional} or \cite[Chapter 18, Proposition 8]{hk-book}; again, the arguments extend without difficulty to arbitrary $\Gamma$.  (One can use Proposition \ref{reverse-system}(ii) to first model $\mathrm{X}$ by a concrete Lebesgue $\Gamma$-system before applying these arguments.)
\end{proof}

To proceed further it is convenient to lift an arbitrary $\Gamma$ to a torsion-free group, and also reduce to the case when $K$ is a Lie group.  The key step is then to establish

\begin{theorem}[Conze--Lesigne equation]\label{conze-eq}  Let $\Gamma$ be a torsion-free countable abelian group, $Z$ an ergodic metrizable rotational $\Gamma$-system, and $K$ a compact abelian Lie group.  Let $\rho$ be a $(Z,K)$-cocycle.  Then the following are equivalent:
\begin{itemize}
    \item[(i)]  $\rho$ is of type $2$.
    \item[(ii)]  $\rho$ obeys the Conze--Lesigne equation (see Definition \ref{cocyc-def}(ix)).
\end{itemize}
\end{theorem}

This result was obtained in \cite[Lemma 8.1]{host2005nonconventional} (or \cite[Section 18.3.3]{hk-book}) in the case when $\Gamma=\Z$ and $K$ is a torus (a connected compact abelian Lie group), building upon previous results in this direction in \cite{cl3}, \cite{fw}.  These arguments extend without much difficulty to arbitrary  torsion-free $\Gamma$ in the connected case when $K$ is a torus; however the case of disconnected $K$ requires additional arguments (cf., the remark after \cite[Lemma C.5]{host2005nonconventional} and Remark \ref{tdK}). The crucial additional case to consider is that of a cyclic group $K = \frac{1}{N} \Z/\Z$.  Here we can proceed instead by some linearization arguments based on those in \cite{btz}.  We remark that thanks to an example of Rudolph \cite{rudolph}, the above theorem fails if the compact abelian Lie group $K$ is replaced with other non-Lie groups, such as solenoid groups; see Remark \ref{rudolph-ex}.

The proof of Theorem \ref{conze-eq} in the general case will be given in Section \ref{conze-sec}.  The derivation of Theorem \ref{main-thm} once Theorem \ref{conze-eq} is in hand is fairly standard and is given in Section \ref{mainproof-sec}, though there are some subtleties requiring the introduction of a topological model in order to properly define the notion of a stabilizer of a certain transitive group action (see \cite{hk-errata} for further discussion of this point).  Here we will use a topological model (which we call a \emph{Koopman model}) that is constructed using Gelfand duality (and the Riesz representation theorem), without the need to impose any ``countability'' conditions such as separability; see Appendix \ref{koopman-sec}.

In order to finish up by expressing an arbitrary Conze--Lesigne system as an inverse limit of nilpotent translational systems, one needs a technical result (Proposition \ref{good-directed}) which states, roughly speaking, that the class of nilpotent translational systems in Theorem \ref{main-thm} is closed under joinings.  As it turns out, this can be established without too much difficulty by exploiting both directions of the equivalences established in Theorem \ref{cl-extend} and Theorem \ref{conze-eq}; see Section \ref{join-sec}.

\subsection{Towards a second-order Pontryagin duality?}\label{sec-speculation}

In principle, Theorem \ref{main-thm} provides a complete description of all ergodic Conze--Lesigne $\Gamma$-systems associated to a given countable abelian group $\Gamma$.  However, one could seek a more tractable such description, in which every  Conze--Lesigne system is described by certain algebraic data from which one can easily answer questions about such systems, such as whether two such systems are isomorphic (or whether one is a factor of the other), whether the system is a translational system or a nilsystem, whether it is generated by a cocycle obeying the Conze--Lesigne equation, whether the structure groups $Z, K$ are connected, and so forth.  

In the case of  Kronecker $\Gamma$-systems, these questions can all be readily answered through Pontryagin duality.  Given a Kronecker $\Gamma$-system $Z$, one can associate the group $E$ of eigenvalues of the system, that is to say those homomorphisms $c \in \hat \Gamma$ from $\Gamma$ to $\T$ such that one has a non-trivial function $f \in L^2(Z)$ for which
$$ T^\gamma f = e(c(\gamma)) f$$
almost everywhere for all $\gamma \in \Gamma$.  This group $E$ is then a subgroup of $\hat \Gamma$ (and is also isomorphic as a group to $\hat Z$).  Conversely, given any subgroup $E$ of $\hat \Gamma$, one can form an associated Kronecker system $Z_E$, defined to be the closure in the compact group $\T^E$ of the subgroup $\phi(\Gamma)$, where $\phi \colon \Gamma \to \T^E$ is the homomorphism
$$ \phi(\gamma) \coloneqq (c(\gamma))_{c \in E},$$
with the rotational $\Gamma$-action given by $\phi$; one can show that this is a Kronecker $\Gamma$-system with eigenvalue group $E$.  The arguments used to establish Theorem \ref{kron} can be used to show that two Kronecker $\Gamma$-systems $Z,Z'$ are isomorphic if and only if their corresponding subgroups $E,E'$ of $\hat \Gamma$ agree (and more generally, $Z$ is a factor of $Z'$ if $E$ is a subgroup of $E'$), thus giving a complete description of Kronecker $\Gamma$-systems in terms of subgroups of $\hat \Gamma$; indeed the above constructions produce a duality of categories.  Furthermore, other properties of the Kronecker $\Gamma$-system can be translated into properties of these subgroups by the usual dictionary of Pontryagin duality.  For instance, given a Kronecker $\Gamma$-system $Z$ and its associated subgroup $E \leq \hat \Gamma$:
\begin{itemize}
    \item $Z$ is separable (or metrizable) if and only if $E$ is countable.
    \item $Z$ is connected if and only if $E$ is torsion-free.
    \item $Z$ is a Lie group if and only if $E$ is finitely generated.
\end{itemize}
In analogy with this state of affairs in order one, one could hope for a ``second order Pontryagin duality'' in which one could associate to each Conze--Lesigne $\Gamma$-system $\mathrm{X}$ some algebraic data (analogous to the subgroup $E$) which determines the isomorphism class of $\mathrm{X}$, as well as other properties of the system, such as whether it is a translational system or a nilsystem, or whether it is associated to a cocycle that obeys the Conze--Lesigne equation.  Ideally, all existing results about such systems (including Theorem \ref{main-thm}) could then be re-interpreted as specific facets of this duality.   We do not at present have a formal proposal for such a duality, though it seems plausible that some form of group cohomology will be involved. We hope to investigate these issues further in subsequent work.

\subsection{Acknowledgments}

AJ was supported by DFG-research fellowship JA 2512/3-1 and by T\"UBITAK-project 123F122. OS was supported by ERC grant ErgComNum 682150, and ISF grant 2112/20. TT was supported by a Simons Investigator grant, the James and Carol Collins Chair, the Mathematical Analysis \& Application Research Fund Endowment, and by NSF grant DMS-1764034.  We thank Bryna Kra and Bernard Host for several discussions of key points in \cite{hk-book}, leading in particular to the errata \cite{hk-errata}.

\section{Notation}\label{notation-sec}

We use $\T \coloneqq \R/\Z$ to denote the additive unit circle.  Given any locally compact abelian group $G = (G,+)$, we define the Pontryagin dual $\hat G$ to be the collection of all continuous homomorphisms from $G$ to $\T$; as is well known, this is also a locally compact abelian group (with the compact-open topology) with $\hat{\hat G} \equiv G$.  We let $S^1$ be the unit circle in $\C$, and let $e\colon \T \to S^1$ be the standard character $e(\theta) \coloneqq e^{2\pi i \theta}$.

We briefly recall the construction of the Host--Kra--Ziegler factors from \cite[\S 3]{host2005nonconventional}, \cite[Chapter 9.1]{hk-book}, or \cite[Appendix A]{btz}.  We begin with the traditional setting of concrete Lebesgue $\Gamma$-systems, with $\Gamma$ a countable abelian group.  Given such a $\Gamma$-system $\mathrm{X} = (X, \mathcal{X}, \mu, T)$, we can recursively define the Host--Kra parallelepiped $\Gamma$-systems $\mathrm{X}^{[k]} = (X^{[k]}, \mathcal{X}^{[k]}, \mu^{[k]}, T^{[k]})$ for $k \geq 0$ by setting $\mathrm{X}^{[0]} \coloneqq \mathrm{X}$ and
$$ \mathrm{X}^{[k+1]} \coloneqq \mathrm{X}^{[k]} \times_{\mathrm{Z}^0(\mathrm{X}^{[k]})} \mathrm{X}^{[k]}$$
where the right-hand side is the relatively independent product of $\mathrm{X}^{[k]}$ with itself over the invariant factor $\mathrm{Z}^0(\mathrm{X}^{[k]})$; see \cite[Chapter 5]{furstenberg2014recurrence} for the construction of relatively independent product for concrete Lebesgue spaces.  As a set, $\mathrm{X}^{[k]}$ can be viewed as a subset of $X^{\{0,1\}^k}$, which we can split as $X \times X^{\{0,1\}^k\backslash \{0\}^k} $.  We then define the Host--Kra--Ziegler factor $\mathrm{Z}^k(\mathrm{X})$ (up to equivalence) by declaring a set $A \in \X$ to be measurable with respect to the $\sigma$-algebra of $\mathrm{Z}^k(\mathrm{X})$ if and only if there is a measurable subset $B$ of $X^{\{0,1\}^k \backslash \{0\}^k}$ such that
\begin{equation}\label{1ab}
 1_A( x_0 ) = 1_B(x_*)
 \end{equation}
for $\mu^{[k]}$-almost all $(x_0,x_*) \in X^{[k]}$.  We refer the reader to \cite[\S 3]{host2005nonconventional}, \cite[Chapter 9]{hk-book}, or \cite[Appendix A]{btz} for the basic properties of these factors, and in particular for their relationship with the Host--Kra--Gowers seminorms (which we will not utilize here).  

Exactly the same constructions can be performed in the more general setting of abstract $\Gamma$-systems, with no requirement of separability (with the analogue of \eqref{1ab} being that the abstract indicator functions $1_A, 1_B$ agree when they are both pulled back to $X^{[k]}$, which is now a probability algebra); alternatively, one can use the canonical model (see Proposition \ref{reverse-system}(iii)) to model such an abstract $\Gamma$-system by a compact $\Gamma$-system and repeat the previous constructions without significant modification.  Note from \cite[Theorem 8.1]{jt-foundational} that the relatively independent product construction is also valid in this ``uncountable'' setting.  


\subsection{Cocycles and extensions}\label{cocyc}

We recall some standard notations for (measurable, abelian) cocycles (largely following \cite{hk-book}, but extended to arbitrary countable abelian groups $\Gamma$).  Here it is convenient to work in the category of concrete $\Gamma$-systems, but permit the cocycles to only be defined up to almost everywhere equivalence (so that the extension generated by such a cocycle is merely an abstract $\Gamma$-system rather than a concrete one).

\begin{definition}[Cocycles and extensions]\label{cocyc-def}  Let $\Gamma$ be a countable abelian group, let $\mathrm{Y} = (Y, {\mathcal Y}, \nu, S)$ be a (concrete) ergodic Lebesgue $\Gamma$-system, and let $K = (K,+)$ be a compact abelian group written additively.
\begin{itemize}
    \item[(i)]  A \emph{$(\mathrm{Y}, K)$-cocycle} is a collection $(\rho_\gamma)_{\gamma \in \Gamma}$ of (concrete) measurable maps $\rho_\gamma \colon Y \to K$ (defined up to almost everywhere equivalence) obeying the \emph{cocycle equation}
    \begin{equation}\label{cocycle-eq}
    \rho_{\gamma_1 + \gamma_2} = \rho_{\gamma_1} \circ S^{\gamma_2} + \rho_{\gamma_2}
    \end{equation}
    $\nu$-almost everywhere for all $\gamma_1,\gamma_2 \in \Gamma$.  Observe that the space of $(\mathrm{Y}, K)$-cocycles form an abelian group under pointwise addition.
    \item[(ii)]  Let $M(\mathrm{Y},K)$ denote the collection of measurable maps from $\mathrm{Y}$ to $K$, up to equivalence almost everywhere; we give this space the topology of convergence in measure (and also endow this space with the Borel $\sigma$-algebra).  If $F \in M(\mathrm{Y},K)$, we define the \emph{derivative} $dF = ((dF)_\gamma)_{\gamma \in \Gamma}$ to be the $(\mathrm{Y}, K)$-cocycle
    $$ (dF)_\gamma \coloneqq F \circ S^\gamma - F.$$
    It is easy to verify that this is indeed a $(\mathrm{Y}, K)$-cocycle.  Any $(\mathrm{Y}, K)$-cocycle of the form $dF$ will be called a \emph{$(\mathrm{Y}, K)$-coboundary}.  Two $(\mathrm{Y}, K)$-cocycles $\rho, \rho'$ are said to be \emph{$(\mathrm{Y}, K)$-cohomologous} if they differ by a $(\mathrm{Y}, K)$-coboundary with respect to the group structure on the space of $(\mathrm{Y}, K)$-cocycles, in which case we write
    $$ \rho \sim_{\mathrm{Y},K} \rho'.$$
    Thus for instance $\rho$ is a $(\mathrm{Y}, K)$-coboundary if and only if $\rho$ is $(\mathrm{Y}, K)$-cohomologous to zero:  $\rho \sim_{\mathrm{Y},K} 0$.
    \item[(iii)]  If $\rho$ is a $(\mathrm{Y},K)$-cocycle, we define the \emph{abelian extension} $\mathrm{Y} \rtimes_\rho K$ to be the concrete probability space that is the product of $(Y, {\mathcal Y}, \nu)$ and $K$ (with the latter equipped with the Haar probability measure), with a \emph{near-action} $T$ given by
    $$ T^\gamma (y, k) \coloneqq (S^\gamma y, k + \rho_\gamma(y))$$
    for all $\gamma \in \Gamma$, $y \in Y$ and $k \in K$, where for each $\gamma$ we arbitrarily select one representative $\rho_\gamma \colon Y \to K$ of the equivalence class for this cocycle.  While from Fubini's theorem one easily sees that each $T^\gamma$ is measure-preserving, the homomorphism law $T^{\gamma_1+\gamma_2} = T^{\gamma_1} \circ T^{\gamma_2}$ for $\gamma_1,\gamma_2 \in \Gamma$ is only true \emph{almost everywhere} rather than everywhere.  Thus, $\mathrm{Y} \rtimes_\rho K$ is not quite well-defined as a concrete $\Gamma$-system; however, as discussed in Appendix \ref{dynamic}, it defines an \emph{abstract} $\Gamma$-system without difficulty (and this system does not depend on the choice of representative of each $\rho_\gamma$).  We say that the $(\mathrm{Y},K)$-cocycle is \emph{ergodic} if this abstract $\Gamma$-system $\mathrm{Y} \rtimes_\rho K$ is ergodic.
    \item[(iv)]  Let $\rho$ be a $(\mathrm{Y},K)$-cocycle.  If $\phi \colon K \to K'$ is a continuous homomorphism from $K$ to another compact abelian group $K'$, we let $\phi \circ \rho$ be the $(\mathrm{Y}, K')$-cocycle
    $$ (\phi \circ \rho)_\gamma \coloneqq \phi \circ \rho_\gamma;$$
    one easily verifies that this is indeed a $(\mathrm{Y}, K')$-cocycle.  Similarly, if $\pi \colon \mathrm{Y}' \to \mathrm{Y}$ is a (concrete) factor map, we let $\rho \circ \pi$ be the $(\mathrm{Y}', K)$-cocycle
    $$ (\rho \circ \pi)_\gamma \coloneqq \rho_\gamma \circ \pi,$$
    which one again easily verifies to be a $(\mathrm{Y}',K)$-cocycle.
    \item[(v)]  If $\rho$ is a $(\mathrm{Y},K)$-cocycle and $V$ is an automorphism of the concrete $\Gamma$-system $\mathrm{Y}$ (thus $V \colon Y \to Y$ is a measure-preserving invertible map and $V \circ S^\gamma = S^\gamma \circ V$ for all $\gamma \in \Gamma$), we define the derivative $\partial_V \rho$ to be the $(\mathrm{Y},K)$-cocycle 
    $$ \partial_V \rho \coloneqq \rho \circ V - \rho.$$
    \item[(vi)]  We let $\Hom(\Gamma,K)$ be the collection of all homomorphisms $c \colon \Gamma \to K$. Every homomorphism $c \in \Hom(\Gamma,K)$ can be viewed\footnote{In the notation of (iv), we are identifying $c$ with $c \circ \mathrm{pt}$, where $\mathrm{pt}$ is the factor map from $\mathrm{Y}$ to a point.} as a $(\mathrm{Y},K)$-cocycle by the formula
    $$ c_\gamma(y) \coloneqq c(\gamma)$$
    for all $\gamma \in \Gamma$ and $y \in Y$. 
    \item[(vii)]  A $(\mathrm{Y},K)$-cocycle $\rho$ is a \emph{$(\mathrm{Y},K)$-quasi-coboundary} if it is  $(\mathrm{Y},K)$-cohomologous to a homomorphism, that is to say there exists measurable $F \colon \mathrm{Y} \to K$ and a homomorphism $c \colon \Gamma \to K$ such that
    $$ \rho_\gamma(y) = F(S^\gamma y) - F(y) + c(\gamma)$$
    for all $\gamma \in \Gamma$ and $\nu$-almost every $y \in Y$.
    \item[(viii)]  If $\rho$ is a $(\mathrm{Y},K)$-cocycle and $k \geq 0$ is an integer, we let $\Delta^{[k]} \rho$ be the $(\mathrm{Y}^{[k]},K)$-cocycle
    $$ (\Delta^{[k]} \rho)_{\gamma}( (y_\omega)_{\omega \in \{0,1\}^k} ) \coloneqq \sum_{\omega \in \{0,1\}^k} (-1)^{|\omega|} \rho_\gamma( y_\omega )$$
    where $|(\omega_1,\dots,\omega_k)| \coloneqq \omega_1 + \dots + \omega_k$.  One easily verifies that $\Delta^{[k]} \rho$ is a $(\mathrm{Y}^{[k]},K)$-cocycle.  If $\Delta^{[k]} \rho$ is a $(\mathrm{Y}^{[k]},K)$-coboundary, we say that $\rho$ \emph{is of type (at most) $k$}.
    \item[(ix)]  If $\mathrm{Y} = Z$ is a rotational system and $\rho$ is a $(Z,K)$-cocycle, we say that $\rho$ \emph{obeys the Conze--Lesigne equation} if for every $z \in Z$, the derivative $\partial_{V_z} \rho$ is a quasi-coboundary, where $V_z$ denotes the translation action $V_z \colon z' \mapsto z+z'$ on $Z$.  In other words, for every $z \in Z$ there exists a measurable $F_z \colon \mathrm{Y} \to K$ and a homomorphism $c_z \colon \Gamma \to K$ such that
    $$ \rho_\gamma(z+z') - \rho_\gamma(z') = F_z(S^\gamma z') - F_z(z') + c_z(\gamma)$$
    for all $\gamma \in \Gamma$ and $\mu_Z$-almost every $z' \in Z$.
\end{itemize}
If the group $K = (K,\cdot)$ is written multiplicatively instead of additively, we define all the preceding concepts analogously, changing all additive notation to multiplicative notation as appropriate.
\end{definition}
In Section \ref{examples-sec}, we collect several examples of Conze--Lesigne systems. In particular these include some concrete examples of cocycles of type $2$ satisfying the Conze--Lesigne equation.  
\begin{remark}
The fact that measurable cocycles only generate abstract $\Gamma$-systems rather than concrete ones will cause some technical issues for us later in our arguments, but these will be resolved by the introduction of suitable topological models, loosely following \cite[\S 19.3.1]{hk-book}, \cite{hk-errata}.
\end{remark}

We recall some basic properties of cocycles:

\begin{proposition}[Basic properties of cocycles]\label{basic-cocycle} Let $\Gamma$ be a countable abelian group, let $\mathrm{Y} = (Y, {\mathcal Y}, \nu, S)$ be an (concrete) ergodic Lebesgue $\Gamma$-system, let $K = (K,+)$ be a compact abelian group, and let $\rho$ be a $(\mathrm{Y},K)$-cocycle.
\begin{itemize}
    \item[(i)] (Moore--Schmidt theorem)  We have $\rho \sim_{\mathrm{Y},K} 0$ if and only if $\xi \circ \rho \sim_{\mathrm{Y},\T} 0$ for all $\xi \in \hat K$.
    \item[(ii)]  (Criterion for ergodicity)  $\rho$ is ergodic if and only if $\xi \circ \rho \not \sim_{\mathrm{Y},\T} 0$ for all non-zero $\xi \in \hat K$.
    \item[(iii)] (Mackey--Zimmer theorem)  If $\mu$ is a $\Gamma$-invariant ergodic probability measure on $Y \times K$ that pushes down to $\nu$ on $Y$, then there exists a closed subgroup $H$ of $K$ (called the \emph{Mackey group} of $\rho, \nu$) and an ergodic $(\mathrm{Y},H)$-cocycle $\rho'$ such that $\rho' \sim_{\mathrm{Y},K} \rho$, and that the  $\Gamma$-system $Y \times_\rho K$ equipped with the measure $\mu$ is abstractly isomorphic to $Y \rtimes_{\rho'} H$ (equipped with product measure).
    \item[(iv)]  (Shifting to be ergodic) If $\Gamma$ is torsion-free and $K$ is connected metrizable, then there exists $c \in \Hom(\Gamma,K)$ such that $\rho + c$ is ergodic. 
    \item[(v)]  (Differentiation lowers type) If $\mathrm{Y} = Z$ is a rotational $\Gamma$-system, $K = \T$, and $\rho$ is of type $2$, then $\partial_{V_z} \rho$ is of type $1$ for all $z \in Z$.
    \item[(vi)]  (Order $1$ cocycles and quasi-coboundaries)  If $\Gamma$ is torsion-free and $K = \T$, then $\rho$ is of type $1$ if and only if $\rho$ is a $(\mathrm{Y},\T)$-quasi-coboundary.
\end{itemize}
\end{proposition}
\begin{proof}  For (i), see \cite[Theorem 4.3]{moore1980coboundaries} or \cite[Theorem 1.1]{jt19} (see also \cite[Chapter 5, Lemma 7]{hk-book} for the $\Gamma=\Z$ case).  For (ii), see \cite[Chapter 5, Lemma 8]{hk-book} (this is stated for $\Gamma=\Z$, but the proof extends without difficulty to arbitrary countable abelian $\Gamma$).    For (iii), see \cite{mackey}, \cite[Corollary 3.8, Theorem 4.3]{zimmer}, \cite[Theorem 3.26]{glasner2015ergodic}, or \cite[Theorem 1.6]{jt20}.  We remark that (iii) is closely related to (i), (ii); for instance, $H$ is the annihilator of the group of characters $\xi \in \hat K$ for which $\xi \circ \rho \sim_{\mathrm{Y},\T} 0$.

Part (iv) is a routine generalization of \cite[Chapter 5, Corollary 9]{hk-book}; for the convenience of the reader we review the argument here.  As $K$ is connected metrizable, the Pontryagin dual $\hat K$ is countable and torsion-free, while $\Hom(\Gamma,K)$ is a compact abelian group.  By (ii), it thus suffices to show that for every non-zero $\xi \in \hat K$, one has $\xi \circ (\rho + c)\not \sim_{\mathrm{Y},\T} 0$ for almost all $c \in \Hom(\Gamma,K)$.  Fixing $\xi$, it suffices upon subtraction to show that $\xi \circ c \not \sim_{\mathrm{Y},\T} 0$ for almost all $c \in \Hom(\Gamma,K)$.  But a character $\xi \circ c \in \hat \Gamma$ is a $(\mathrm{Y},\T)$-coboundary if and only if it is an eigenvalue of the action $S$; from the separability of $\mathrm{Y}$, there are countably many such eigenvalues, so it suffices to show that $\xi \circ c \neq 0$ for almost all $c \in \Hom(\Gamma,K)$.  If this is not the case, then the closed subgroup $\{ c \in \Hom(\Gamma,K): \xi \circ c = 0 \}$ would have finite index in $\Hom(\Gamma,K)$, hence there is an integer $n$ such that $\xi \circ n c = 0$ for all $c \in \Hom(\Gamma,K)$.  As $\Gamma$ is torsion-free, this would imply that $n\xi$ vanishes, contradicting the torsion-free nature of $\hat K$.

Part (v) was established for $\Z$-actions in \cite[Corollary 7.5(i)]{host2005nonconventional} or \cite[Chapter 18, Proposition 11(i)]{hk-book}, and for general actions\footnote{There is a typo in the statement of that lemma: the hypothesis that $X$ be of order $<k$ should instead be that $Z_{<k}(X)$ be a factor of $Y$.} in \cite[Lemma 5.3]{btz} by the same method.

Part (vi) was established for $\Z$-actions in \cite[Chapter 5, Lemma 13]{hk-book}, but the extension to torsion-free $\Gamma$ is routine: for the convenience of the reader we review the proof here.  The ``if'' part is easy, so we focus on the ``only if'' part.  By using (iv) to shift $\rho$ by a character (which does not affect $\Delta^{[1]} \rho$) we may assume without loss of generality that $\rho$ is ergodic.  By hypothesis, $\Delta^{[1]} \rho$ is a $(\Y^{[1]}, \T)$-coboundary, thus there exists a measurable map $F \colon Y \times Y \to \T$ obeying the equation
$$ \Delta^{[1]} \rho_\gamma  = F \circ S^\gamma - F$$
for all $\gamma \in \Gamma$.  Setting $\mathrm{X} \coloneqq \mathrm{Y} \rtimes_\rho \T$, we conclude that the map
$$ H \colon ((y_0,k_0), (y_1,k_1)) \mapsto  e(F(y_0,y_1) + k_1-k_0)$$
is $\Gamma$-invariant in $\mathrm{X} \times \mathrm{X}$.  Thus the integral operator $T_H$ with kernel $H$ is a non-trivial Hilbert-Schmidt operator on $L^2(\mathrm{X})$ that commutes with the $\Gamma$-action, thus there is an eigenfunction $\phi \in L^2(\mathrm{X})$ of this action that is also a non-trivial eigenfunction of $T_H^* T_H$.  The function $\beta \in L^2(\mathrm{Y})$ defined by
$$ \beta(y) \coloneqq \int_\T \bar \phi(y,k) e(k)\ dk$$
cannot vanish identically (since otherwise $T_H \phi$ would vanish), and obeys the equation
$$ \beta(S^\gamma y) = e(-c(\gamma)) e(\rho_\gamma(y)) \beta(y)$$
for all $\gamma \in \Gamma$ and almost every $y$, where $e(c(\gamma))$ is the eigenvalue of $\phi$ with respect to $S^\gamma$; note that $c$ is necessarily a character in $\hat \Gamma$.  The function $|\beta|$ is $\Gamma$-invariant, hence constant.  Writing $\beta = |\beta| e(\tilde F)$, we obtain $\rho_\gamma = \tilde F \circ S^\gamma - \tilde F + c(\gamma)$, thus $\rho_\gamma$ is a $(\mathrm{Y},\T)$-quasi-coboundary as desired.
\end{proof}

\begin{remark}\label{tdK} As remarked  after \cite[Lemma C.5]{host2005nonconventional}, Proposition \ref{basic-cocycle}(vi) can fail if the circle $\T$ is replaced by a disconnected group.  For instance, take $\Gamma = \Z$, $K = \frac{1}{2}\Z/\Z$, and $\mathrm{Y}$ to be the rotational ergodic separable $\Z$-system $\T$ with action $\phi(n) = n \alpha \hbox{ mod } 1$ for all $n \in \Z$ and some irrational real $\alpha$.  If we let $\{ \} \colon \T \to [0,1)$ be the fractional part map, then one can check that the tuple $\rho = (\rho_n)_{n \in \Z}$ given by
$$ \rho_n(x) \coloneqq \frac{ \{ x+n\alpha\} - \{x\} - n \alpha }{2} \hbox{ mod } 1$$
is a $(\Z,K)$-cocycle.  From the identity
$$ \Delta^{[1]} \rho = d F^{[1]}$$
where $F^{[1]} \in M(\T^2, K)$ is the function
$$ F^{[1]}(x,y) \coloneqq \frac{\{x\}-\{y\} - \{x-y\}}{2} \hbox{ mod } 1$$
we see that the $(\Z,K)$-cocycle is of order $1$; however it is not a $(\Z,K)$-quasi-coboundary.  Indeed, if there was some $F \in M(\T,K)$ and $c \in \Hom(\Z,K)$ such that $\rho = dF + c$, then by specializing to $n=1$ we conclude that
$$  \frac{ \{ x+\alpha\} - \{x\} - \alpha }{2} = 
F(x+\alpha)-F(x) + c(1) \hbox{ mod } 1$$
or equivalently
$$ f(x+\alpha) = e\left(\frac{\alpha}{2}+c(1)\right) f(x)$$
for all $x \in \T$, where 
$$f(x) \coloneqq e\left( F(x)-\frac{\{x\}}{2}\right).$$
By Fourier analysis this implies that $\frac{\alpha}{2}+c(1)$ needs to be an integer multiple of $\alpha$, but this is inconsistent with the irrationality of $\alpha$ since $c(1) \in \frac{1}{2} \Z/\Z$.

On the other hand, if the system $\mathrm{Y}$ is \emph{$2$-divisible} in the sense that its Kronecker factor has a divisible Pontryagin dual, then one can replace the circle $\T$ in Proposition \ref{basic-cocycle}(vi) by an arbitrary compact abelian group $K$; this follows from that proposition and \cite[Proposition 3.8]{shalom2}.  As a consequence, the requirement in Theorem \ref{conze-eq} that $K$ be a Lie group can be dropped in the $2$-divisible case; this is essentially \cite[Theorem 4.1]{shalom2}.
\end{remark}

\section{Derivation of the Conze--Lesigne equation}\label{conze-sec}

In this section we establish Theorem \ref{conze-eq}.  

We begin with the derivation of (i) from (ii).  Observe from Definition \ref{cocyc-def}(ix) that if the $(Z,K)$-cocycle $\rho$ obeys the Conze--Lesigne equation, then the $(Z,\T)$-cocycle $\xi \circ \rho$ obeys the Conze--Lesigne equation for any $\xi \in \hat K$.  Similarly, from the Moore--Schmidt theorem (Proposition \ref{basic-cocycle}(i)), Definition \ref{cocyc-def}(viii) and the obvious identity $\Delta^{[k]}(\xi \circ \rho) = \xi \circ \Delta^{[k]}(\rho)$ for any $\xi \in \hat K$ and $k \geq 0$ we see that if $\xi \circ \rho$ is of type $k$ for every $\xi \in \hat K$, then $\rho$ is of type $k$.  From these observations we see that to show that (ii) implies (i), it suffices to do so in the case $K=\T$. 
By Definition \ref{cocyc-def}(ix), we see that for every $z \in Z$ there exists a $F_z \in M(Z,\T)$ and a character $c_z \in \hat \Gamma$ such that
\begin{equation}\label{rho-eq}
 \partial_{V_z} \rho = d F_z + c_z.
 \end{equation}
At this point we run into the technical issue that $F_z$ and $c_z$ need not depend in a measurable fashion on $z$.  It is however possible to select $F_z, c_z$ so that this is the case, by means of the following result:

\begin{proposition}[Measurable selection]\label{mes-select}  Let $\Gamma$ be a countable abelian group, let $\mathrm{Y}$ be a concrete ergodic Lebesgue $\Gamma$-system, and let $U$ be a measurable space.  Suppose that we have a measurable map $u \mapsto h_u$ from $U$ to the space of $(\mathrm{Y},\T)$-cocycles (which we can view as a subset of $M(\mathrm{Y},\T)^\Gamma$, which we endow with the product topology), with the property that for each $u \in U$ we can find $F_u \in M(\mathrm{Y},\T)$ and a character $c_u \in \hat \Gamma$ such that
$$ h_u = d F_u + c_u.$$
Then, after adjusting $F_u$ and $c_u$ as necessary, we may ensure that $F_u, c_u$ depend in a measurable fashion on $u$.
\end{proposition}

This proposition is a special case of \cite[Lemma C.4]{btz} (which handles a more general situation in which the $h_u$ need not obey the cocycle equation, and the $c_u$ are allowed to be polynomials of a given degree).  The proof of that lemma requires at one point the measurability of a certain function $n_u$ constructed in that proof.  The verification of this measurability is actually somewhat non-trivial, and so we give a complete proof of Proposition \ref{mes-select} in Appendix \ref{mes-select-app}.  As remarked in \cite{btz}, there are several other ways to establish this proposition, including using Borel cross sections of homomorphisms between Polish groups (see \cite[Theorem A.1]{host2005nonconventional}) or a general measurable selection result of Dixmier (cf. \cite[Theorem 1.2.4]{kechrisbecker}). In the case $\Gamma=\Z$ this result was essentially established in \cite[Proposition 10.5]{fw}.

Invoking Proposition \ref{mes-select}, we can now select the $F_z, c_z$ solving \eqref{rho-eq} to depend in a measurable fashion on $z$.  

As observed in \cite[\S 3.2]{host2005nonconventional} (see also \cite[\S 8.1.2]{hk-book}), $Z^{[2]}$ can be viewed as a translational system on the compact group
$$Z^{[2]} = \{ (z, z+s_1,z+s_2,z+s_1+s_2): z,s_1,s_2 \in Z\}$$
with translation action $\phi^{[2]} \colon \Gamma \to Z^{[2]}$ given by the diagonal action
$$ \phi^{[2]}(\gamma) \coloneqq (\phi(\gamma), \phi(\gamma), \phi(\gamma), \phi(\gamma))$$
and $\phi \colon \Gamma \to Z$ the original translation action on $Z$.  The $(Z^{[2]},\T)$-cocycle $\Delta^{[2]} \rho$ is then given by the formula
\begin{align*}
 (\Delta^{[2]} \rho)_\gamma(z, z+s_1,z+s_2,z+s_1+s_2) &= \rho_\gamma(z) - \rho_\gamma(z+s_1) - \rho(z+s_2) + \rho(z+s_1+s_2) \\
 &= \partial_{V_{s_2}} \rho_\gamma(z) -  \partial_{V_{s_2}} \rho_\gamma(z+s_1).
\end{align*}
Applying \eqref{rho-eq}, we conclude the identity
$$ \Delta^{[2]} \rho = d F^{[2]}$$
where $F^{[2]} \colon  Z^{[2]} \to \T$ is the function
$$ F^{[2]}(z,z+s_1, z+s_2, z+s_1+s_2) \coloneqq F_{s_2}(z) - F_{s_2}(z+s_1).$$
By construction of the $F_z$, $F^{[2]}$ is measurable, and hence by Definition \ref{cocyc-def}(viii) $\rho$ is of type $2$.  This concludes the derivation of (i) from (ii).

Now we show that (i) implies (ii).  Any compact abelian Lie group is isomorphic to the direct product of a torus and a finite abelian group (see e.g., \cite[Exercise 1.4.27(iii)]{tao-hilbert}), and hence also isomorphic to the direct product of finitely many copies of the circle $\T$ and cyclic groups $\frac{1}{N}\Z/\Z$.  It is clear that to show (i) implies (ii) for a direct product $K = K_1 \times K_2$, it suffices to do so for the two factors $K_1$ and $K_2$ separately.  Thus it suffices to establish this implication in the special cases $K = \T$ and $K = \frac{1}{N}\Z/\Z$ for a natural number\footnote{Using the Chinese remainder theorem one could reduce further to the case when $N$ is a power of a prime, but this does not seem to simplify the argument significantly.} $N$.

The $K=\T$ case is immediate from Proposition \ref{basic-cocycle}: if $z \in Z$, then Lemma \ref{basic-cocycle}(v) implies that $\partial_{V_z} \rho$ is of type $1$, and Lemma \ref{basic-cocycle}(vi) then gives that $\partial_{V_z} \rho$ is a $(Z,\T)$-quasi-coboundary, thus giving the required Conze--Lesigne equation.

We turn to the $K = \frac{1}{N}\Z / \Z$ case.  Now one cannot directly apply Lemma \ref{basic-cocycle}(v), (vi).  However, since $K$ is a subgroup of $\T$, we can also view the $(Z,K)$-cocycle $\rho$ as a $(Z,\T)$-cocycle, which will of course still be of type $2$.  Applying the previous argument, we conclude that $\partial_{V_z} \rho$ is a $(Z,\T)$-quasi-coboundary for every $z \in \T$, thus we can find $c_z \in \Hom(\Gamma,\T) = \hat \Gamma$ such that
\begin{equation}\label{caz}
\partial_{V_z} \rho \sim_{Z,\T} c_z.
\end{equation}
By Proposition \ref{mes-select}, we may ensure that $c_z$ depends in a measurable fashion on $z$.

The main difficulty here is that the homomorphism $c_z$ takes values in $\T$ rather than in the smaller group $K$.  To resolve this, we need some additional structural control on the $c_z$.  We first apply a translation $V_{z'}$ to \eqref{caz} to conclude that
$$ (\partial_{V_z} \rho) \circ V_{z'} \sim_{Z,\T} c_z$$
for any $z,z' \in Z$; combining these identities with the cocycle identity
$$\partial_{V_{z+z'}} \rho = (\partial_{V_z} \rho) \circ V_{z'} + \partial_{V_{z'}} \rho$$
we conclude that
$$ c_{z+z'} - c_z - c_{z'} \sim_{Z,\T} 0.$$
Thus, if we let $E \leq \hat \Gamma$ denote the group
$$ E \coloneqq \{ c \in \hat \Gamma: c \sim_{Z,\T} 0 \}$$
then the map $z \mapsto c_z$ is a homomorphism mod $E$, in the sense that
\begin{equation}\label{czz}
 c_{z+z'} = c_z + c_{z'} \hbox{ mod } E
 \end{equation}
for all $z,z' \in Z$.

Note that if $c \in E$, then $c = dF$ for some $F \in M(Z,\T)$, which implies that $e(F)$ is an eigenfunction of the rotational system $Z$:
$$ e(F) \circ T^\gamma = e(c(\gamma)) e(F).$$
By the unitary nature of the action, eigenfunctions with different eigenvalues are orthogonal.  Since $L^2(Z)$ is separable, we conclude that $E$ is countable.

We can now locally remove the ``mod $E$'' reduction in \eqref{czz} by the following argument (cf. the proof of \cite[Proposition 6.1]{btz}).  By \eqref{czz}, the map $(z,z') \mapsto c_{z+z'}-c_z-c_{z'}$ is a measurable map from $Z \times Z$ to the countable set $E$.  The autocorrelation function
$$ a \mapsto \mu_{Z^2}( \{ (z,z') \in Z^2: c_{z+z'}-c_z-c_{z'} = c_{z+a+z'}-c_{z+a}-c_{z'} \}),$$
where $\mu_{Z^2}$ is the Haar probability measure on $Z^2$, is then a continuous function on $Z$ which equals $1$ at $0$ (this follows for instance from Lusin's theorem).  Thus there exists an open neighborhood $U$ of the identity such that
$$ \mu_{Z^2}( \{ (z,z') \in Z^2: c_{z+z'}-c_z-c_{z'} = c_{z+a+z'}-c_{z+a}-c_{z'} \}) \geq 0.9$$
(say) for all $a \in U$.  Canceling the $c_{z'}$ and making the change of variables $z'' = z+z'$, we see that for all $a \in U$, we have
$$ c_{z+a} - c_z = c_{z''+a}- c_{z''}$$
for at least $0.9$ of pairs $(z,z'') \in Z^2$ by measure, which implies that there exists a (unique) $c'_a \in \hat \Gamma$ such that
\begin{equation}\label{cza}
 c_{z+a} - c_z = c'_a
 \end{equation}
for at least $0.9$ of the $z \in Z$ by measure; furthermore, $c'_a$ will depend measurably on $a$ (it is the mode of $c_{z+a}-c_z$).  From \eqref{czz} we see that
\begin{equation}\label{caa}
 c'_a = c_a \hbox{ mod } E
 \end{equation}
for all $a \in U$, and from several applications of \eqref{cza} we have
\begin{equation}\label{cia}
c'_a + c'_b = c'_{a+b}
\end{equation}
whenever $a,b,a+b \in U$.

We return to the equation \eqref{caz}.  Since $\rho$ takes values in $K = \frac{1}{N}\Z / \Z$, we have $N \rho = 0$, hence from \eqref{caz} $Nc_z \sim_{Z,\T} 0$ for all $z \in Z$, hence by \eqref{caa} we have $Nc'_a \in E$ for all $a \in U$.  Thus there is $e \in E$ such that $Nc'_a = e$ for all $a$ in a positive measure subset of $U$; from \eqref{cia} and the Steinhaus lemma, we conclude that $Nc'_a = 0$ for all $a$ in an open neighborhood $U' \subset U$
of the identity.  Thus for $a \in U'$, $c'_a$ takes values in $K$, and so from \eqref{caz}, \eqref{caa} $\partial_{V_a} \rho - c'_a$ is a $(Z,K)$-cocycle which is a $(Z,\T)$-coboundary.  By the Moore--Schmidt theorem (Proposition \ref{basic-cocycle}(i)), $\partial_{V_a} \rho - c'_a$ is also a $(Z,K)$-coboundary (note that all the characters of $K$ are of the form $k \mapsto nk$ for some integer $n$).  Thus $\partial_{V_a} \rho$ is a $(Z,K)$-quasi-coboundary for all $a \in U'$.  Meanwhile, from the identity
$$ \partial_{V_{\phi(\gamma)}} \rho = d \rho_{\gamma}$$
for any $\gamma \in \Gamma$, we see that $\partial_{V_a} \rho$ is also a $(Z,K)$-quasi-coboundary (in fact a $(Z,K)$-coboundary) for all $a \in \phi(\Gamma)$.  By the cocycle identity
$$ \partial_{V_{a+b}} \rho = (\partial_{V_a} \rho) \circ V_b + \partial_{V_b} \rho$$
for any $a,b \in Z$, we conclude that $\partial_{V_a} \rho$ is a $(Z,K)$-quasi-coboundary for all $a \in \phi(\Gamma) + U'$.  But since the rotational $\Gamma$-system $Z$ is ergodic, the subgroup $\phi(\Gamma)$ of $Z$ is dense, and hence $\phi(\Gamma) + U'$ is all of $Z$.  Thus $\partial_{V_a} \rho$ is a $(Z,K)$-quasi-coboundary for all $a \in Z$, thus $\rho$ obeys the Conze--Lesigne equation.  This concludes the derivation of (ii) from (i), and the proof of Theorem \ref{conze-eq} is complete.

\section{Conclusion of the argument}\label{mainproof-sec}

In this section we establish Theorem \ref{main-thm}.

\subsection{From nilpotent translational systems to the Conze--Lesigne equation}\label{nil-cl}

We begin with the derivation of (i) from (ii).  From Lemma \ref{basic-facts}(iii), Theorem \ref{cl-extend} and Theorem \ref{conze-eq}, it suffices to show the following claim:

\begin{proposition}[Verifying the Conze--Lesigne equation]\label{verify}  
Let $\Gamma$ be a countable abelian group, and let $G/\Lambda$ be an ergodic translational $\Gamma$-system, where $G$ is a locally compact nilpotent Polish group of nilpotency class $2$, and $\Lambda$ is a lattice in $G$, and one also has a closed central subgroup $G_2$ of $G$ containing $[G,G]$ such that $\Lambda \cap G_2$ is a lattice in $G_2$.   Then $G/\Lambda$ is abstractly isomorphic to a group extension $Z \rtimes_\rho K$, where $Z$ is a rotational $\Gamma$-system, $K$ is a compact abelian group, and $\rho$ is a $(Z, K)$-cocycle obeying the Conze--Lesigne equation.  (Note that $Z$ is not required to be the Kronecker factor.)
\end{proposition}

We now prove this proposition. We let $\phi \colon \Gamma \to G$ denote the translation action. We take $Z$ to be the compact group $G/G_2 \Lambda$, written additively.  We write $\pi \colon G \to Z$ for the projection homomorphism; this map factors through the quotient map from $G$ to $G/\Lambda$, and we use $\tilde \pi \colon G/\Lambda \to Z$ to denote the projection map produced in this fashion.  Then $Z$ is a rotational $\Gamma$-system with action given by $\pi \circ \phi$.  Next, we take $K$ to be the compact group $G_2 / (G_2 \cap \Lambda)$ written additively.  Because $G_2$ is central, this group $K$ acts freely on $G/\Lambda$; we express this action additively, thus if $k \in K$ and $x \in G/\Lambda$, we write $k+x=x+k$ for the action of $k$ on $x$.  By construction, we thus have
\begin{equation}\label{action}
g_2 x = x + \Pi(g_2)
\end{equation}
whenever $x \in G/\Lambda$ and $g_2 \in G_2$, where $\Pi \colon G_2 \to K$ is the projection homomorphism.  Observe that the orbits of this free $K$-action on $G/\Lambda$ are precisely the fibers of $\tilde \pi$, thus $G/\Lambda$ is a principal $K$-bundle over $Z$ (as a set, at least). Also, by the central nature of $G_2$ we see that 
\begin{equation}\label{gxk}
g(x+k) = gx + k
\end{equation}
 for all $g \in G$, $x \in G/\Lambda$, and $k \in K$.

It will be convenient to ``work in coordinates'' to facilitate computations.
We claim that the projection map $\tilde \pi \colon G/\Lambda \to Z$ admits a Borel cross-section, that is to say a Borel-measurable map $s \colon Z \to G/\Lambda$ such that $\tilde \pi(s(z)) = z$ for all $z \in Z$.  Indeed, observe that the map $\pi \colon G \to Z$ is a continuous surjective homomorphism of Polish groups, hence by \cite[Theorem 1.2.4]{kechrisbecker} this map admits a Borel cross-section $s' \colon Z \to G$; quotienting out by $\Lambda$ then gives the claim.

For any $\gamma \in \Gamma$ and $z \in Z$, the points $\phi(\gamma) s(z)$ and $s(\pi \circ \phi(\gamma) + z)$ in $G/\Lambda$ both lie in the fiber $\tilde \pi^{-1}(\pi \circ \phi(\gamma) + z)$, so there is a unique element $\rho_\gamma(z)$ of $K$ for which one has the identity
\begin{equation}\label{phis}
 \phi(\gamma) s(z) = s(\pi \circ \phi(\gamma) + z) + \rho_\gamma(z).
\end{equation}
It is easy to see that $\rho_\gamma \colon Z \to K$ is measurable for each $\gamma$.  By computing $\phi(\gamma_1+\gamma_2) s(z)$ in two different ways using \eqref{gxk},\eqref{phis} we can verify that $\rho \coloneqq (\rho_\gamma)_{\gamma \in \Gamma}$ is in fact a $(Z,K)$-cocycle.  By the identification
$$ (z,k) \equiv s(z) + k$$
one can then verify that the translational system $G/\Lambda$ is abstractly isomorphic to the semi-direct product $Z \rtimes_\rho K$ (one can use Fubini's theorem to check that the product measure of $Z \rtimes_\rho K$ is invariant under the left action of $G$ under this identification and is thus identified with the Haar probability measure on $G/\Lambda$).  

To conclude the proof of Proposition \ref{verify}, it will suffice to show that $\rho$ obeys the Conze--Lesigne equation, which can be achieved by standard calculations in a suitable coordinate system as follows.  Let $z_0 \in Z$ be arbitrary.  As the projection $\pi \colon G \to Z$ is surjective, we can find\footnote{For this argument we will not need to require $g_{z_0}$ to depend in a measurable fashion on $z_0$, though we could ensure this if desired, by using either a variant of the Borel section constructed previously, or by a variant of Proposition \ref{mes-select}.} $g_{z_0} \in G$ such that $\pi(g_{z_0}) = z_0$.  Applying $g_{z_0}$ to \eqref{phis} and using \eqref{gxk}, we see that
$$ g_{z_0} \phi(\gamma) s(z) = g_{z_0} s(\pi \circ \phi(\gamma) + z) + \rho_\gamma(z)$$
for any $\gamma \in \Gamma$ and $z \in Z$.  Writing $g_{z_0} \phi(\gamma) = [g_{z_0},\phi(\gamma)] \phi(\gamma) g_{z_0}$ and noting that the commutator $[g_{z_0},\phi(\gamma)]$ lies in $G_2$, we then have from \eqref{action} that
\begin{equation}\label{phig}
 \phi(\gamma) g_{z_0} s(z) + \Pi( [g_{z_0},\phi(\gamma)] ) = g_{z_0} s(\pi \circ \phi(\gamma) + z) + \rho_\gamma(z).
 \end{equation}
On the other hand, since $g_{z_0} s(z)$ and $s(z+z_0)$ both lie in the fibre $\tilde \pi^{-1}(z+z_0)$, there exists a unique measurable function $F_{z_0} \colon Z \to K$ such that
$$ g_{z_0} s(z) = s(z+z_0) + F_{z_0}(z)$$
for all $z \in Z$.  Inserting this (both for $z$ and for $\pi \circ \phi(\gamma) + z$) into equation \eqref{phig} and using \eqref{gxk} and \eqref{phis}, we conclude that
$$ s(\pi \circ \phi(\gamma) + z + z_0) + \rho_\gamma(z+z_0) + F_{z_0}(z) + \Pi( [g_{z_0},\phi(\gamma)] ) = 
 s(\pi \circ \phi(\gamma) + z + z_0) + F_{z_0}(\pi \circ \phi(\gamma) + z) + \rho_\gamma(z);$$
 as the $K$-action is free, this can be rearranged as
 $$\rho_\gamma(z+z_0) - \rho_\gamma(z) = F_{z_0}(\pi \circ \phi(\gamma) + z)  -  F_{z_0}(z)  - \Pi( [g_{z_0},\phi(\gamma)] )$$
 or equivalently
\begin{equation}\label{cs}
 \partial_{V_{z_0}} \rho = d F_{z_0} + c_{z_0}
\end{equation}
where $c_{z_0} \colon \Gamma \to K$ is the map
$$ c_{z_0}(\gamma) \coloneqq - \Pi( [g_{z_0},\phi(\gamma)] ).$$
As $G$ has nilpotency class $2$, one easily verifies that $c_{z_0}$ is a homomorphism (this also follows from the fact that the other terms in \eqref{cs} are $(Z,K)$-cocycles).  Hence $\rho$ obeys the Conze--Lesigne equation as required.  This completes the proof of Proposition \ref{verify}, and hence the derivation of (i) from (ii).

\begin{remark}\label{g2} For a given translational $\Gamma$-system $G/\Lambda$, with $G$ a nilpotent locally compact Polish group of nilpotency class two, there can be some flexibility in how to select the subgroup $G_2$; it must contain the commutator group $[G,G]$ and be contained in turn in the center $Z(G)$ of $G$, and needs to be closed and ``rational'' in the sense that $G_2 \cap \Lambda$ is cocompact in $G_2$, but is otherwise arbitrary.  From the above discussion, this means that it is possible for the translational $\Gamma$-system $G/\Lambda$ to be expressed as an abelian group extension $Z \rtimes_\rho K$ of a rotational $\Gamma$-system by a cocycle $\rho$ obeying the Conze--Lesigne equation in several inequivalent ways.  As discussed in Remark \ref{kron-remark} below, the minimal choice $G_2 = [G,G]$ corresponds to the case when $Z$ is the maximal rotational $\Gamma$-system factor of $G/\Lambda$, i.e., the Kronecker factor; however in some cases one can also take larger choices of $G_2$, such as the center $Z(G)$ of $G$, which correspond to smaller choices of $Z$.  See Section \ref{first-ex} for one example of this situation.
\end{remark}

\subsection{From the Conze--Lesigne equation to nilpotent translational systems}\label{cl-nil}

Now we show that (i) implies (ii).  We begin with a technical reduction.  Define a \emph{good system} to be a translational $\Gamma$-system $G/\Lambda$ of the form required in part (ii) with the additional conditions listed at the end of the theorem; thus $G$ is a locally compact nilpotent Polish group of nilpotency class two, $\Lambda$ is an abelian lattice in $G$, and $G$ contains a compact central Lie group $G_2$ containing $[G,G]$ with $\Lambda \cap G_2$ trivial.  Call a factor of a $\Gamma$-system a \emph{good factor} if it is abstractly isomorphic to a good system.  Our task is to show that any Conze--Lesigne $\Gamma$-system is the inverse limit of a directed family of good factors.  The requirement to be a directed set can be dropped thanks to the following observation (cf., \cite[\S 13.3.2, Proposition 16]{hk-book}):

\begin{proposition}[Good factors form a directed set]\label{good-directed}  Given two good factors $\mathrm{Y}_1, \mathrm{Y}_2$ of $\mathrm{X}$, there exists another good factor $\mathrm{Y}$ of $\mathrm{X}$ such that $\mathrm{Y}_1, \mathrm{Y}_2  \leq \mathrm{Y}$.
\end{proposition}

We defer the proof of this proposition to Section \ref{join-sec}.  Assuming it for now, any family of good factors can be completed to a directed set of good factors by a transfinite induction on the cardinality of the family of good factors. Indeed, the base case being trivial, suppose $\beta=\alpha+1$ is a successor ordinal and $\mathfrak{X}_\alpha$ is a directed family of good factors. Applying Proposition \ref{good-directed}, we can form $\mathfrak{X}_\beta$ by taking the join of each good factor in $\mathfrak{X}_\alpha$ with the additional element. Moreover, if $\beta$ is a limit ordinal, then we can take the union of all $\mathfrak{X}_\alpha$ over $\alpha<\beta$. 

Hence for the purposes of showing that (i) implies (ii) we can now drop the requirement that the family of good factors be directed.  In particular, if $\mathrm{X}$ is the inverse limit of some other systems $\mathrm{X}_n$, and each $\mathrm{X}_n$ was already demonstrated to be an inverse limit of good factors, then $\mathrm{X}$ itself must also be an inverse limit of good factors, simply by concatenating all the families of good factors together (and ignoring the directed set requirement).

The next step is to reduce to the case of separable $\Gamma$-systems.  By Lemma \ref{basic-facts} and the preceding discussion, it suffices to show that every $\Gamma$-system $\mathrm{X}$ is the inverse limit of separable $\Gamma$-systems.  But given any finite collection ${\mathcal F}$ of elements of the $\sigma$-complete Boolean algebra $\X$ associated to $\mathrm{X}$, one can form the factor $\mathrm{X}_{\mathcal F}$ by replacing $\X$ with the $\sigma$-complete subalgebra generated by the $\Gamma$-orbit $\{ T^\gamma F: F \in {\mathcal F}, \gamma \in \Gamma \}$ of ${\mathcal F}$, and restricting the measure and action appropriately.  It is clear that this is a separable factor of $\mathrm{X}$, and $\mathrm{X}$ is the inverse limit of the $\mathrm{X}_{\mathcal F}$, as claimed.

Henceforth $\mathrm{X}$ is separable.
The group $\Gamma$ is not assumed to be torsion-free, but it is of course isomorphic to a quotient $\Gamma'/\Sigma$ of a torsion-free countable abelian group $\Gamma'$.  For instance one can take $\Gamma' = \bigoplus_{\gamma \in \Gamma} \Z$ to be the free abelian group formally generated by the elements $\gamma$ of $\Gamma$, with $\Gamma$ then naturally identified with the quotient of $\Gamma'$ by the subgroup $\Sigma$ consisting of formal integer combinations of elements of $\Gamma$ that sum to zero.  Any ergodic separable $\Gamma$-system $\mathrm{X}$ can then be viewed as an ergodic $\Gamma'$-system in the obvious fashion; and if $\mathrm{X}$ when viewed as a $\Gamma'$-system is the inverse limit of good $\Gamma'$-systems $G_n/\Lambda_n$, then each of the factor $\Gamma'$-systems $G_n/\Lambda_n$ must have a trivial action of $\Sigma$ and thus also be interpretable as a good $\Gamma$-system.  Furthermore, if $\mathrm{X}$ is of order $k$ as a $\Gamma$-system for a given $k$, it is easy to see from the definitions that it is also of order $k$ when viewed as a $\Gamma'$-system.  Thus, to prove the implication of (ii) from (i) for $\Gamma$, it suffices to do so for $\Gamma'$.  In particular, we may now assume without loss of generality that $\Gamma$ is torsion-free.

By Theorem \ref{cl-extend}, we can assume without loss of generality that the ergodic separable $\Gamma$-system $\mathrm{X}$ is of the form $\mathrm{X} = Z \rtimes_\rho K$, where $Z$ is a rotational ergodic $\Gamma$-system, $K$ is a compact abelian group, and $\rho$ is a $(Z,K)$-cocycle of type $2$.  Since $\mathrm{X}$ was separable, $Z$ is also separable, hence by Pontryagin duality $\hat Z$ is countable and $Z$ is metrizable.  By the Peter--Weyl theorem or Pontryagin duality (see e.g., \cite[Exercise 1.4.26]{tao-hilbert}), $K$ is the inverse limit of compact abelian Lie groups $K_n$.  One can then easily verify that $Z \rtimes_\rho K$ is the inverse limit of $Z \rtimes_{\pi_n \circ \rho} K_n$, where $\pi_n \colon K \to K_n$ are the projection homomorphisms.  Since $\rho$ is a $(Z,K)$-cocycle of type $2$, the $(Z,K_n)$-cocycles $\pi_n \circ \rho$ also have type $2$.  Thus, for the purposes of establishing Theorem \ref{main-thm}(i), we may assume without loss of generality that $K$ is a compact abelian Lie group.  In particular, by Theorem \ref{conze-eq}, the $(Z,K)$-cocycle $\rho$ now obeys the Conze--Lesigne equation. Also, by Pontryagin duality, $\hat K$ is a finitely generated discrete group.

To summarize so far, we have reduced the derivation of (ii) from (i) to establishing the following proposition (which can be viewed as a partial converse to Proposition \ref{verify}): 
\begin{proposition}[Constructing a nilpotent translational system]\label{construct}  Let $\Gamma$ be a countable abelian group, $Z$ a metrizable rotational ergodic $\Gamma$-system, $K$ a compact abelian Lie group, and $\rho$ an ergodic $(Z,K)$-cocycle obeying the Conze--Lesigne equation.  Then the (ergodic, separable) $\Gamma$-system $Z \rtimes_\rho K$ is abstractly isomorphic to a good system, i.e., a translational system $G/\Lambda$ with $G$ is a locally compact nilpotent Polish group of nilpotency class two, $\Lambda$ is an abelian lattice in $G$, and $G$ contains a compact central Lie group $G_2$ containing $[G,G]$ with $\Lambda \cap G_2$ trivial.
\end{proposition}

To avoid circularity in our arguments we emphasize that our proof of Proposition \ref{construct} will not use Proposition \ref{good-directed}, as this latter proposition has not yet been proven.

To prove Proposition \ref{construct}, we now follow a standard construction (see  \cite{cl2}, \cite{cl3}, \cite{meiri}, \cite{rudolph}, \cite{hk1}, \cite{hk2}, \cite{host2005nonconventional}, \cite{ziegler2007universal}), but taking care to keep track of which structures are only defined abstractly (or up to almost everywhere equivalence), rather than pointwise.
Define the \emph{Host--Kra group} $G$ to be the collection of pairs $(u,F)$, where $u \in Z$ and $F \in M(Z,K)$ obeys the Conze--Lesigne equation
\begin{equation}\label{cle}
\partial_{V_u} \rho = dF + c
\end{equation}
for some homomorphism $c \colon \Gamma \to K$.  These pairs $(u,F)$ generate a near-action  on $\mathrm{X} = Z \rtimes_\rho K$ by the formula
\begin{equation}\label{ufz}
(u,F) (z,k) \coloneqq (z+u, k+F(z)),
\end{equation}
where we arbitrarily select one concrete representative $F \colon Z \to K$ from the equivalence class of $F$.   
One verifies from Fubini's theorem that this near-action is concretely measure-preserving, and that the abstract action on $\mathrm{X}$ does not depend on the choice of representative.  
Thus if we endow $G$ with the group law
$$ (u,F) (u',F') \coloneqq (u+u', F \circ V_{u'} + F')$$
and inverse operation
$$ (u,F)^{-1} \coloneqq (-u, -F \circ V_{-u})$$
one easily verifies that $G$ is a group that has a near-action on $\mathrm{X}$, and thus (as discussed in Appendix \ref{dynamic}) acts \emph{abstractly} on $\mathrm{X}$.  We claim that this abstract action is faithful.  Indeed, if $(u,F)$ acts abstractly trivially on $\mathrm{X}$, then for every bounded measurable  $f \colon Z \times K \to \R$ we have
$$ f( z+u, k+F(z) ) = f( z, k)$$
for $Z \times K$-almost every $(z,k)$. Testing this against functions of the form $f(z,k) = e(\chi(z))$ for characters $\chi \in \hat Z$, we conclude that $u$ vanishes; testing against functions of the form $f(z,k) = e(\xi(k))$ for characters $\xi \in \hat K$ we conclude that $F$ vanishes almost everywhere, giving the claim.  Thus we can identify $G$ with a subgroup of the unitary group on $L^2(\mathrm{X})$, by identifying each $(u,F) \in G$ with the Koopman operator defined in $L^2(\mathrm{X})$ by the usual formula
$$ ((u,F) f)(z,k) \coloneqq f( (u,F)^{-1} (z,k) )$$
for $f \in L^2(\mathrm{X})$)
(note that this is well-defined as a unitary map on $L^2(\mathrm{X})$ that does not depend on the choice of representatives for $F$ or $f$).  

By identifying $k\in K$ with the constant function $z\mapsto k$, we see that $(0,k)$ obeys the Conze--Lesigne equation \eqref{cle}, and hence the group
$$ G_2 \coloneqq \{ (0,k): k \in K \}$$ is a central subgroup of $G$.  We also claim that $G_2$ contains $[G,G]$.  Indeed, if $(u,F), (u', F') \in G$, then a brief calculation shows that
\begin{equation}\label{brief}
[(u,F), (u',F')] = (0, \tilde F)
\end{equation}
where $\tilde F \coloneqq (\partial_{V_{u'}} F - \partial_{V_u} F') \circ V_{-u-u'}$.  Differentiating the formula for $\tilde F$ using \eqref{cle} we see that $d\tilde F = 0$, and hence by ergodicity $\tilde F = k$ for some $k \in K$, giving the required inclusion.  In particular, $G$ is nilpotent with nilpotency class at most two.  From the Conze--Lesigne equation we see that the projection map $(u,F) \mapsto u$ is a surjective homomorphism from $G$ to $Z$ with kernel $H \coloneqq \{ (u,F) \in G: u = 0\}$ (which is clearly an abelian group containing $G_2$ as a subgroup), thus we have the short exact sequence
\begin{equation}\label{short-exact}
0 \to H \to G \to Z \to 0.
\end{equation}
If we let $\phi \colon \Gamma \to G$ denote the map
$$ \phi(\gamma) \coloneqq (\phi_Z(\gamma), \rho_\gamma)$$
where $\phi_Z \colon \Gamma \to Z$ is the rotation action on $Z$, one checks from the definitions that the abstract action of $\Gamma$ on $\mathrm{X}$ is the composition of $\phi$ with the abstract action of $G$ on $\mathrm{X}$, thus
$$ T^\gamma = \phi(\gamma)$$
as abstract maps on $\mathrm{X}$.

 The strong operator topology gives the structure of a Hausdorff topological group to the group of unitary operators in $L^2(\mathrm{X})$, and hence also to $G$.  This is a good topological structure to place on $G$:

\begin{proposition}  $G$ is a locally compact Polish group, and $G_2$ is a closed subgroup of $G$ (and thus also locally compact Polish).
\end{proposition}

\begin{proof} As $\mathrm{X}$ is a separable probability algebra, the Hilbert space $L^2(\mathrm{X})$ is also separable. As is well known, the group of unitary operators on such a space, when equipped with the strong operator topology, is a Polish group.  To show that $G$ is also a Polish group, it thus suffices to show that $G$ is closed in the strong operator topology.  But if $(u_n,F_n) \in G$ is a Cauchy sequence in the strong operator topology, it is easy to see (by testing against characters $\chi \in \hat \Z$) that $u_n$ is a Cauchy sequence in $Z$ that must therefore converge to some $u \in Z$, and for any character $\xi \in \hat K$, that the $\xi \circ F_n$ are a Cauchy sequence in measure, so (by the finitely generated nature of $\hat K$) $F_n$ converges in measure to some limit $F \colon Z \to K$.  It is then not difficult to show that $(u,F)$ obeys the Conze--Lesigne equation and that $(u_n,F_n)$ converges to $(u,F)$, which demonstrates that $G$ is closed and thus a Polish group. The same argument shows that $G_2$, $H$ are closed subgroups of $G$; as the obvious bijection from $K$ to $G_2$ is a continuous map, we conclude that $G_2$ is isomorphic to $K$ as a compact abelian group.  In particular $G, G_2$ are second countable.

It remains to show that $G$ is locally compact.  The homomorphism from $G$ to $Z$ is a continuous surjective homomorphism of Polish groups, and is thus an open map (see e.g., \cite{kechrisbecker}).  Using the short exact sequence \eqref{short-exact}, we conclude that $G/H$ is isomorphic to $Z$ and is in particular locally compact.  To show that $G$ is locally compact, it thus suffices (see \cite[Theorem 5.25]{hewittross}) to show that $H$ is locally compact.

Since $G_2$ is already compact, it suffices to show that $G_2$ is an open subgroup of $H$, or equivalently that every sequence $(0,F_n)$ in $H$ converging to the identity lies in $G_2$ for $n$ large enough.  By the Conze--Lesigne equation \eqref{cle}, the $F_n$ obey the equation 
\begin{equation}\label{dfc}
dF_n + c_n = 0
\end{equation}
for some $c_n \in \Hom(Z,K)$, so for each $\xi \in \hat K$ we have
\begin{equation}\label{xifn}
 e(\xi(F_n)) \circ S^\gamma = e(\xi(c_n(\gamma)) e(\xi(F_n)) 
 \end{equation}
almost everywhere for all $n$ and all $\gamma \in \Gamma$, where $S$ denotes the rotation action on $Z$.  By the preceding discussion, $F_n$ converges in measure to zero, so for any fixed $\xi \in \hat K$, $e(\xi(F_n))$ converges in measure to $1$.  In particular, for $n$ large enough, $e(\xi(F_n))$ has mean one.  Integrating \eqref{xifn}, we conclude that $\xi(c_n(\gamma)) = 0$ for sufficiently large $n$ and all $\gamma \in \Gamma$.  Since $\hat K$ is finitely generated, we conclude that $c_n=0$ for all sufficiently large $n$.  Thus by \eqref{dfc}, for all sufficiently large $n$, $F_n$ is $\Gamma$-invariant, and therefore constant by ergodicity.  In other words, $(0,F_n)$ lies in $G_2$, giving the claim.
\end{proof}

Note that as $G$ is locally compact and nilpotent, it is unimodular.
It remains to show that $\mathrm{X}$ is abstractly isomorphic to a translational $\Gamma$-system $G/\Lambda$ for some lattice $\Lambda$ in $G$, with $\Lambda \cap G_2$ a lattice in $G_2$.  If $G$ acted concretely (or better yet, continuously) on $\mathrm{X}$, one could hope to proceed here by showing that the action of $G$ on $\mathrm{X}$ was transitive, and take $\Lambda$ to be the stabilizer of a point.  Unfortunately, the action of $G$ that we have on $\mathrm{X}$ is only an abstract action.  To resolve this we use the Koopman topological model $\hat{\mathrm{X}} = (\hat X, \cdot)$ of the abstract $G$-system $\mathrm{X}$ constructed in Theorem \ref{koop}, where we will use $g \colon \hat x \mapsto g \hat x$ to denote the $G$-action on this model.  By Lemma \ref{translation}, it now suffices to establish the following claims:

\begin{itemize}
    \item[(iv)] For any $\hat x_1, \hat x_2 \in \hat X$, there exists $g \in G$ such that $g \hat x_1 = \hat x_2$.
    \item[(v)] For some $\hat x_0 \in \hat X$, the stabilizer $\Lambda \coloneqq \{ g \in G: g \hat x_0 = \hat x_0\}$ is a lattice in $G$, $\Lambda \cap G_2$ is trivial (and hence a lattice in $G_2$), with $\Lambda$ abelian.
\end{itemize}

Indeed, Lemma \ref{translation} will then guarantee that $\mathrm{X}$ is abstractly isomorphic as a $G$-system to the translational $G$-system $G/\Lambda$, and then by applying the group homomorphism $\phi \colon \Gamma \to G$ we see that $\mathrm{X}$ and $G/\Lambda$ are abstractly isomorphic as $\Gamma$-systems as well.

We begin with (iv).  Observe from \eqref{ufz} (and the continuity of the projection from $G$ to $Z$) that any continuous function $f \in C(Z)$ on $Z$ pulls back to a $G$-continuous function $(z,k) \mapsto f(z)$ on $\mathrm{X}$, where the notion of $G$-continuity was defined in Theorem \ref{koop}. Thus we have a tracial $C^*$-algebra homomorphism from $C(Z)$ to the algebra ${\mathcal A}$ of $G$-continuous functions, which preserves the $G$-action (letting $(u,F)$ act on $Z$ by translation by $u$).  By Gelfand--Riesz duality (see \cite[Theorem 5.11]{jt-foundational}) and Theorem \ref{koop}, we thus have a continuous factor map $\hat \pi \colon \hat X \to Z$ of compact $G$-systems.  Because the projection of $G$ to $Z$ is surjective, the (continuous) action of $G$ on $Z$ is transitive.  Thus to establish the transitivity property (iv), it suffices to do so in a single fiber of $\hat \pi$, that is to say we may assume without loss of generality that $\hat \pi(\hat x_1) = \hat \pi(\hat x_2)$.

It suffices to establish transitivity of the $G_2$-action on fibers of $\hat \pi$, that is to say that under the hypothesis $\hat \pi(\hat x_1) = \hat \pi(\hat x_2)$ there exists $k \in K$ such that $(0,k) \hat x_1 = \hat x_2$.  We now repeat the arguments from \cite[\S 19.3.3, Lemma 10]{hk-book}.  Suppose for contradiction that the $G_2$-orbit of $\hat x_1$ does not contain $\hat x_2$, then by continuity we can find an open neighborhood $U$ of $\hat x_1$ in $\hat X$ such that $\hat x_2$ does not lie in the $G_2$-orbit $\{ (0,k) \hat x: \hat x \in U; k \in K \}$ of $U$.  By Urysohn's lemma, we can find a non-negative function $f \in C(\hat X)$ supported on $U$ that is positive at $\hat x_1$; the averaged function
$$\overline{f}(\hat x) \coloneqq \int_K f((0,k) \hat x)\ dk,$$
with $dk$ the Haar probability measure on $K$, is then a $G_2$-invariant function in $C(\hat X)$ that is non-zero at $\hat x_1$ but vanishes at $\hat x_2$.  By construction of the Koopman model, $\overline{f}$ can then be identified with a $G$-continuous function in $L^\infty(X)$ which is also $G_2$-invariant, and hence arises from a $G$-continuous function on $Z$ thanks to \eqref{ufz}.  But as the projection from $G$ to $Z$ is surjective, the $G$-continuous functions on $Z$ can be identified with the ordinary continuous functions on $Z$, thus $\overline{f}$ can be identified with an element of $C(Z)$.  But as $\hat x_1, \hat x_2$ lie in the same fiber of $\hat \pi$ we must then have $\overline{f}(\hat x_1) = \overline{f}(\hat x_2)$, giving the required contradiction.  This establishes the transitivity property (iv).

As a corollary of this transitivity and the faithfulness of the $G$ action, we see (cf., \cite[\S 19.3.3, Lemma 11]{hk-book}) that the central $G_2$ action must be free, since if $(0,k) \hat x = \hat x$ for some $k \in K$ and $\hat x \in \hat X$, then by transitivity and centrality the action of $(0,k)$ on $\hat X$ would be trivial, hence $k=0$.  Thus if we let $\Lambda$ be a stabilizer of a point $\hat x_0$ in $G$, then $\Lambda \cap G_2$ is trivial and thus clearly a lattice in the compact group $G_2$.  Since $[\Lambda,\Lambda]$ is contained in both $\Lambda$ and $[G,G] \subset G_2$, it must be trivial, hence $\Lambda$ is abelian.

To complete the verification of (v) we need to show that the stabilizer group $\Lambda$ is a lattice in $G$.  Since the $G$-action on $\hat X$ projects down to the $G$-action on $Z$, the stabilizer group $\Lambda$ must be contained in the kernel $H$ of the projection from $G$ to $Z$.  On the other hand, as $G_2$ is an open subgroup of $H$ and $\Lambda \cap G_2$ is trivial, we conclude that $\Lambda$ is discrete.  Also, by the transitivity of the $G_2$ action on fibers of $\pi$ (which are preserved by the action of $H$) we see that $\Lambda$ must intersect every coset of $G_2$ in $H$.  Thus the quotient $H/\Lambda$ is homeomorphic to $G_2$ and thus compact.  Since $G/H \equiv Z$ is also compact, the projection from $G$ to $Z$ is open, and $G$ is locally compact, $G/\Lambda$ is also compact\footnote{Indeed, the local compactness of $G$, the open nature of the projection, and the compactness of $G/H$ gives an inclusion $G \subset FH$ for some compact $F$, and the compactness of $H/\Lambda$ gives an inclusion $H \subset F'\Lambda$ for some compact $F'$, thus $G \subset FF' \Lambda$ and hence $G/\Lambda$ is compact.}, so $\Lambda$ is a lattice as required.  This concludes the proof of Proposition \ref{construct}, and hence Theorem \ref{main-thm} once we establish Proposition \ref{good-directed}.

\begin{remark}  The above arguments can also establish an isomorphism $H \equiv K \times \Lambda$ of topological groups; we leave the details to the interested reader.
\end{remark}

\begin{remark}\label{kron-remark} In the model case where $Z$ is the Kronecker factor of $Z \rtimes_\rho K$, we can upgrade the inclusion $[G,G] \subset G_2$ in the above construction to $[G,G] = G_2$ (where by $[G,G]$ we denote the closed group generated by the commutators).  We sketch the proof as follows.  Suppose this claim failed, then by Pontryagin duality there exists a character $\xi \in \hat G_2 \equiv \hat K$ that annihilates $[G,G]$.  For any $z \in Z$, let $F_z, c_z$ be a solution to the Conze--Lesigne equation \eqref{rho-eq}.  Then $(z,F_z)$, $\phi(\gamma) = (\phi_Z(\gamma), \rho_\gamma)$ both lie in $G$, and by \eqref{brief} their commutator is $(0,c_z(\gamma)) \in G_2$, thus $c_z(\gamma)$ is annihilated by $\xi$.  Applying $\xi$ to the Conze--Lesigne equation \eqref{rho-eq}, we conclude that $\partial_{V_z}(\xi \circ \rho)$ is a $(Z,\T)$-coboundary for every $z \in Z$.  If we let $\pi \colon K \to K/[G,G]$ be the quotient homomorphism (identifying $[G,G] \leq G_2$ with a subgroup of $K$ in the obvious fashion), we conclude from the Moore--Schmidt theorem (Proposition \ref{basic-cocycle}(i)) that $\partial_{V_z}(\pi \circ \rho)$ is a $(Z, K/[G,G])$-coboundary for every $z$. By a variant of Proposition \ref{mes-select}, this implies that $\pi \circ \rho$ is of type $1$, and hence (by a variant of Theorem \ref{cl-extend}) $Z \rtimes_{\pi \circ \rho} K/[G,G]$ is of order $1$, i.e., a Kronecker system.  Thus $Z$ is not the maximal rotational factor of $Z \rtimes_\rho K$, giving the required contradiction.
\end{remark}

\subsection{Joinings of good systems}\label{join-sec}

Finally, we supply the proof of Proposition \ref{good-directed}.  It suffices to establish the following claim (cf. \cite[\S 11.2.3, Corollary 10]{hk-book}):

\begin{proposition}[Measure classification on good systems]\label{class-good}  Let $G/\Lambda$ be a (possibly non-ergodic) good system.  If $\nu$ is an $\Gamma$-invariant ergodic measure on $G/\Lambda$, then $G/\Lambda$ equipped with $\nu$ is abstractly isomorphic to a good system.
\end{proposition}

Indeed, suppose that an (abstract) $\Gamma$-system $\mathrm{X}$ had two good factors, which we write without loss of generality as $G_1/\Lambda_1$ and $G_2/\Lambda_2$.  The abstract factor maps give pullback maps from $C(G_1/\Lambda_1)$ and $C(G_2/\Lambda_2)$ to $L^\infty(\mathrm{X})$, which by the Stone--Weierstrass theorem gives a pullback map from $C(G_1 \times G_2/\Lambda_1 \times \Lambda_2)$ to $L^\infty(\mathrm{X})$ which one can verify to be a $C^*$-homomorphism.  The integral on $\mathrm{X}$ then induces a trace on $C(G_1 \times G_2/\Lambda_1 \times \Lambda_2)$, which by the Riesz representation theorem gives a measure $\nu$ on $G_1 \times G_2/\Lambda_1 \times \Lambda_2$ (in fact it gives a joining of $G_1/\Lambda_1$ and $G_2/\Lambda_2$).  By construction, 
$G_1 \times G_2/\Lambda_1 \times \Lambda_2$ equipped with $\nu$ is an ergodic $\Gamma$-system that is a factor of $\mathrm{X}$, and has $G_1/\Lambda_1$ and $G_2/\Lambda_2$ as factors in turn.  By Proposition \ref{class-good}, this factor is a good factor, giving Proposition \ref{good-directed}.

It remains to establish Proposition \ref{class-good}.  This turns out to be a straightforward consequence of the implications regarding Conze--Lesigne systems that we have already established. By lifting $\Gamma$ to a torsion-free group as before, we may assume without loss of generality that $\Gamma$ is torsion-free. By Proposition \ref{verify}, the translational $\Gamma$-system $G/\Lambda$ is abstractly isomorphic to a group extension $Z \rtimes_\rho K$ for some (possibly non-ergodic) separable rotational system $Z$, some compact abelian Lie group $K$, and some cocycle $\rho$ obeying the Conze--Lesigne equation, except that the measure $\nu$ is not necessarily equal to the product measure on $Z \rtimes_\rho K$.  An inspection of the construction shows that the cocycle equation \eqref{cocycle-eq} holds everywhere (not just almost everywhere), and similarly for the Conze--Lesigne equation.  Thus this group extension $Z \rtimes_\rho K$ (which by abuse of notation we also equip with the measure $\nu$) is well-defined as a concrete $\Gamma$-system, not just an abstract one.

The ergodic measure $\nu$ on $Z \rtimes_\rho K$ pushes down to an ergodic measure $\nu_Z$ on $Z$, which is invariant under a translational action $\phi_Z$ of $\Gamma$ on $Z$.  A standard Fourier-analytic computation then shows that $\nu_Z$ must be Haar measure of a coset of some closed subgroup $Z'$ of $Z$ (indeed, $Z'$ is the closure of $\phi_Z(\Gamma)$ in $Z$).  Applying a translation, we may assume without loss of generality that the coset of $Z'$ is just $Z'$ itself.  The ergodic measure $\nu$ is then supported on a subsystem $Z' \rtimes_{\rho'} K$ of $Z \rtimes_\rho K$, where $\rho'$ is the restriction of $\rho$ to $Z'$.  Since $\rho$ obeys the cocycle and Conze--Lesigne equations everywhere (not just almost everywhere), the same is true for $\rho'$; that is to say, $\rho'$ is a $(Z',K)$-cocycle that obeys the Conze--Lesigne equation. By Theorem \ref{cl-extend}, the $(Z',K)$-cocycle $\rho'$ is of type 2.

By construction, the ergodic measure $\nu$ on $Z' \rtimes_{\rho'} K$ pushes down to the Haar measure on $Z'$.  By the Mackey--Zimmer theorem (see Lemma \ref{basic-facts}(iii)), there is a closed subgroup $H$ of $K$ and an ergodic $(Z',H)$-cocycle $\rho''$ such that $\rho''$ is $(Z',K)$-cohomologous to $\rho'$, and $Z' \rtimes_{\rho'} K$ equipped with the measure $\nu$ is abstractly isomorphic to $Z' \rtimes_{\rho''} H$ equipped with product measure.  Since the $(Z',K)$-cocycle $\rho'$ is of type $2$, the $(Z',K)$-cohomologous cocycle $\rho''$ is of type $2$ when viewed as a $(Z',K)$-cocycle, thus $\Delta^{[2]} \rho''$ is a $((Z')^{[2]}, K)$-coboundary.  As 
$\Delta^{[2]} \rho''$ is also a $((Z')^{[2]}, H)$-cocycle, we see from the Moore--Schmidt theorem (Proposition \ref{basic-cocycle}(i)) that $\Delta^{[2]} \rho''$ is a $((Z')^{[2]}, H)$-coboundary (note from Pontryagin duality that every character on $H$ extends (not necessarily uniquely) to a character on $K$).  Thus $\rho''$ is also of type $2$ when viewed as a $(Z',H)$-cocycle.

Since $K$ is a compact abelian Lie group, the closed subgroup $H$ is also a compact abelian Lie group.  Applying Theorem \ref{cl-extend} again, we see that $\rho''$ obeys the Conze--Lesigne equation.  Applying Proposition \ref{construct} (which did not require the use of Proposition \ref{good-directed} in its proof), we conclude that $Z' \rtimes_{\rho''} H$ is a good system.  Since this system is isomorphic to $G/\Lambda$ equipped with the measure $\nu$, the claim follows.

\begin{remark} It is worth considering whether results such as Proposition \ref{good-directed} or Proposition \ref{class-good} can be established directly from the theory of nilpotent translational systems, without relying on the implications presented in Theorem \ref{cl-extend}, Proposition \ref{verify}, or Proposition \ref{construct}. Although \cite[\S 13.3.2, Proposition 16]{hk-book} accomplished this in the case of nilsystems, the arguments presented there heavily rely on the finite dimensionality of these systems. However, a somewhat similar situation arises also for $\Z$-actions, as we currently lack a direct method for proving that a factor of an inverse limit of nilsystems is itself an inverse limit of nilsystems, without relying on the Host--Kra--Ziegler structure theorem.
\end{remark}

\section{Some examples of Conze--Lesigne systems}\label{examples-sec}

In this section we give some examples of Conze--Lesigne systems (in even, odd, and zero characteristic respectively) to illustrate the main theorems.

\subsection{First example: an extension of a characteristic two rotational system}\label{first-ex}

Let $\Gamma \coloneqq \F_2^\omega$ be the countably generated vector space over $\F_2$, and let $Z \coloneqq \F_2^\N$ be the countable product of $\F_2$ equipped with Haar probability measure $\nu=\mu_{\F_2^\N}$, and let $S:\Gamma\to \Aut(Z,\nu)$ be the $\Gamma$-rotation $S^\gamma(z) \coloneqq z+\gamma$ (using the obvious identification of $\Gamma$ with a subgroup of $Z$). 
By the mean ergodic theorem, the projection of any $f\in L^\infty(Z)$ depending only on finitely many coordinates onto the invariant subspace of $L^2(Z)$ is constant. The span of these functions is dense in $L^2(Z)$. 
Hence $(Z,\nu,S)$ is an ergodic separable $\Gamma$-rotational system. 

Let $K \coloneqq \Z/4\Z$ be the cyclic group of order $4$, and let $\rho = (\rho_\gamma)_{\gamma \in \Gamma}$ be the $(Z, K)$-cocycle 
\begin{equation}\label{rhog}
\rho_\gamma(z) \coloneqq \sum_{n \in \N} (-1)^{z_n} 1_{\gamma_n=1},
\end{equation}
where we define $$(-1)^{x}=
\begin{cases}
1, & x=0, \\
-1, & x=1,
\end{cases}$$
and $1_{\gamma_n=1}$ similarly equals $1$ when $\gamma_n=1$ and $0$ otherwise.  It is easy to verify that $\rho$ is a $(Z, K)$-cocycle.  We have the following further properties:

\begin{proposition}[Properties of $\rho$]\label{rho}\
\begin{itemize}
    \item[(i)] $\rho$ is ergodic.
    \item[(ii)]  $\rho$ is of type $2$.
    \item[(iii)]  $\rho$ obeys the Conze--Lesigne equation.
\end{itemize}
\end{proposition}

\begin{proof}
We begin with (i).  By Proposition \ref{basic-cocycle}(ii), it suffices to show that $\xi \circ \rho$ is not a $(Z,\T)$-coboundary for any non-zero $\xi \in \hat K$.  Since the character $x \mapsto \frac{x}{2} \mod 1$ of $K$ is a multiple of any non-zero $\xi \in \hat K$, it suffices to verify the claim for this specific character.  Suppose for contradiction that $\frac{1}{2} \circ \rho \mod 1$ is a coboundary, thus there exists $F \in M(Z,\T)$ such that
\begin{equation}\label{eq-1}
\frac{1}{2} \rho_\gamma = F\circ V_\gamma - F \mod 1
\end{equation}
$\nu$-almost everywhere for all $\gamma$.  In particular we have
$$ F(z + e_n) = F(z) + \frac{1}{2} \mod 1$$
for any generator $e_n$ of $\Gamma = \F_2^\omega$.  But by Lusin's theorem, $F(z+e_n)$ converges in measure to $F(z)$ as $n \to \infty$, giving a contradiction.

Now we verify (ii). 
As observed in \cite[\S 3.2]{host2005nonconventional} (see also \cite[\S 8.1.2]{hk-book}), $Z^{[2]}$ can be viewed as a translational system on the compact group
$$Z^{[2]} = \{ (z, z+s_1,z+s_2,z+s_1+s_2): z,s_1,s_2 \in Z\},$$
with each $\gamma \in \Gamma$ acting by translation by $(\gamma,\gamma,\gamma,\gamma)$.  Thus we need to locate a measurable function $F \colon Z^{[2]} \to \T$ such that
\begin{equation}\label{type-2}
\begin{split}
&    \rho_\gamma(z) - \rho_\gamma(z+s_1) - \rho_\gamma(z+s_2) + \rho_\gamma(z+s_1+s_2)\\
&\quad = F(z+\gamma,z+s_1+\gamma,z+s_2+\gamma,z+s_1+s_2+\gamma) - F(z,z+s_1,z+s_2,z+s_1+s_2) 
\end{split}
\end{equation}
for all $\gamma \in \Gamma$ and $\nu$-almost all $z,s_1,s_2$.  But the left-hand side expands as
$$ \sum_{n \in \N} (-1)^{z_n} (1 - (-1)^{s_{1,n}}) (1 - (-1)^{s_{2,n}}) 1_{\gamma_n=1},$$
and in the group $K = \Z/4\Z$, the product $(1 - (-1)^{s_{1,n}}) (1 - (-1)^{s_{2,n}})$ always vanishes.  Thus we may simply take $F=0$ to verify that $\rho$ is of type $2$.

Finally, we establish (iii).  We need to show that for each $z \in Z$, there exists $F_z \in M(Z,K)$ and $c_z \in \Hom(\Gamma,K)$ such that $\partial_{V_z} \rho = dF_z + c_z$, or in other words that
\begin{equation}\label{cleq}
 \rho_\gamma(w+z) - \rho_\gamma(w) = F_z(w+\gamma) - F_z(w) + c_z(\gamma)
\end{equation}
for all $\gamma \in \Gamma$ and $\nu$-almost all $w \in Z$.  But the left-hand side expands as
$$ \sum_{n \in \N} (-1)^{w_n} ((-1)^{z_n}-1) 1_{\gamma_n = 1},$$
and in the cyclic group $K = \Z/4\Z$, $(-1)^{w_n} ((-1)^{z_n}-1)$ is equal to $(-1)^{z_n}-1$.  Thus we can solve the Conze--Lesigne equation by setting $F_z \coloneqq 0$ and
$$ c_z(\gamma) \coloneqq  \sum_{n \in \N} ((-1)^{z_n}-1) 1_{\gamma_n = 1}$$
which one easily verifies to be a homomorphism from $\Gamma$ to $K$.
\end{proof}

By the above proposition and Theorem \ref{cl-extend}, $Z \rtimes_\rho K$ is an (ergodic, separable) Conze--Lesigne $\Gamma$-system.  Now we compute its Host--Kra group $G$.  By definition, this is the set of all pairs $(u,F)$, where $u \in Z$ and $F \in M(Z,K)$ obeys the Conze--Lesigne equation
$$ \partial_{V_u} \rho = dF + c$$
for some homomorphism $c \in \Hom(\Gamma,K)$.  By the proof of Proposition \ref{rho}(iii), $\partial_{V_u} \rho_\gamma$ is constant, thus $(dF)_\gamma$ is constant for each $\gamma \in \Gamma$.  In particular, for each natural number $n$, there must be a constant $c_n \in K$ such that
\begin{equation}\label{al}
F(z+e_n) - F(z) = c_n
\end{equation}
for almost all $z \in Z$.  Shifting $z$ by $e_n$ and summing in the characteristic two group $Z$, we conclude that $2c_n=0$, thus $c_n$ is either equal to $0$ or $2$.  On the other hand, $F(\cdot+e_n)$ converges in measure to $F$, hence all but finitely many of the $c_n$ vanish.  From this we conclude that $F$ must take the form
\begin{equation}\label{eq-hk-group}F(z) = \theta + \sum_{n \in \N} (-1)^{z_n} 1_{\sigma_n = 1}\end{equation}
almost everywhere for some $\theta \in K$ and $\sigma \in \Gamma$, which are uniquely determined by $F$.  Thus, by abuse of notation, we can write the Host--Kra group $G$ as the collection of triples\footnote{To see the converse that any element of $Z \times K \times \Gamma$ can be identified with an element of $G$, let $(u,\theta,\sigma)\in Z \times K \times \Gamma$. We need to find $c\in \Hom(\Gamma,K)$ such that $\partial_{V_u} \rho = dF + c$ where $F$ is defined by \eqref{eq-hk-group} for the given choice of $(\theta,\sigma)$. For the given $u\in Z$, by Proposition \ref{rho}(iii), there are $F'\in M(Z,K)$ (which we can also represent as in \eqref{eq-hk-group}) and $c' \in \Hom(\Gamma,K)$ such that $\partial_{V_u} \rho = dF' + c'$. By a direct computation, one verifies that $c\coloneqq d(F'-F)+c'\in \Hom(\Gamma,K)$.}  $(u, \theta, \sigma) \in Z \times K \times \Gamma$, and one can calculate the group law in $Z \times K \times \Gamma$ to be
$$ (u,\theta,\sigma) (u',\theta',\sigma') \coloneqq \left(u+u', \theta+\theta' + \sum_{n \in \N} ((-1)^{u'_n}+1) 1_{\sigma_n=1}, \sigma+\sigma'\right)$$
and inverse 
$$ (u,\theta,\sigma)^{-1} = \left(-u, -\theta + \sum_{n \in \N} ((-1)^{u_n}-1) 1_{\sigma_n=1}, -\sigma\right).$$
This group acts transitively (and continuously) on $Z \times K$ by the formula
$$ (u,\theta,\sigma) (z,k) \coloneqq \left(z+u, k + \theta + \sum_{n \in \N} (-1)^{z_n} 1_{\sigma_n = 1}\right)$$
and the stabiliser $\Lambda$ of the point $(0,0)$ is
$$ \Lambda \coloneqq \left\{ \left(0, -\sum_{n \in \N} (-1)^{z_n} 1_{\sigma_n = 1}, \sigma\right): \sigma \in \Gamma \right\}.$$
One can check that the strong operator topology on $G$ corresponds to the product topology on $Z \times K \times \Gamma$ (viewing $K,\Gamma$ as discrete groups), so that $G$ is a second countable locally compact Polish group, $\Lambda$ is a lattice in $G$, and $Z \times K$ is isomorphic (as a compact $\Gamma$-space) to $G/\Lambda$, with the action of a group element $\gamma \in \Gamma$ on $G/\Lambda$ given by multiplication by $(\gamma,0,\gamma)$.  If one defines the subgroup $G_2$ of $G$ by
$$ G_2 \coloneqq \{ (0, \theta, 0): \theta \in K \}$$
then $G_2$ is a closed central subgroup of $G$ that contains\footnote{In fact $[G,G]$ is strictly smaller than $G_2$, consisting only of those triples $(0,\theta,0)$ where $\theta$ is a multiple of two; we leave the verification of this fact to the interested reader.  Compare also with Remark \ref{g2}.} $[G,G]$, and hence $G$ is nilpotent of class two; also $\Lambda \cap G_2$ is trivial and thus a lattice in the compact group $G_2$.   One can now verify that Theorem \ref{main-thm} holds for this example. 

As $Z$ is a rotational system, it is contained in the Kronecker factor of $Z \rtimes_\rho K$, thanks to Theorem \ref{kron}.  However, the Kronecker factor turns out to be slightly larger than this:

\begin{proposition}\label{kroncalc}  The Kronecker factor of $Z \rtimes_\rho K$ is $Z \rtimes_{2\rho} 2K$, where $2K = 2\Z/4\Z$ is a cyclic group of order $2$, with factor map $(z,k) \mapsto (z,2k)$.
\end{proposition}

\begin{proof}  Observe that the action of a group element $\gamma$ on $Z \rtimes_{2\rho} 2K$ is given by translation in the group $Z \times 2K$ by $(\gamma, 2 \sum_{n \in \N} 1_{\gamma_n=1})$.  Thus $Z \rtimes_{2\rho} 2K$ is a translational $\Gamma$-system and thus contained in the Kronecker factor.

To establish the converse claim, observe from Theorem \ref{kron}, Proposition \ref{basic-facts}, and Pontryagin duality that the factor algebra of the Kronecker factor is generated by eigenfunctions of the $\Gamma$-action, that is to say functions $f \in L^\infty(Z \rtimes_\rho K)$ such that
\begin{equation}\label{eigen}
f \circ T^\gamma = \lambda_\gamma f
\end{equation}
almost everywhere for all $\gamma \in \Gamma$ and some $\lambda_\gamma \in \C$ (which must lie in $S^1$ by unitarity).  Since the $\Gamma$-action on $Z \rtimes_\rho K$ commutes with the abelian $K$-action, we see on applying a Fourier decomposition with respect to the $K$ variable that we can restrict attention to eigenfunctions of the form
\begin{equation}\label{fz}
 f(z,k) = F(z) e(mk/4)
 \end{equation}
for some $m=0,1,2,3$ and some $F \in L^\infty(Z)$.  The eigenfunctions with even $m$ already are measurable in the factor $Z \times_\rho K$, so it suffices to show that there are no non-trivial eigenfunctions of the form \eqref{fz} with $m$ odd.  The function $|F|$ is $\Gamma$-invariant, thus constant by ergodicity; we may normalize $|F|=1$.  Applying the eigenfunction equation \eqref{eigen} with $\gamma=e_n$ we see after some calculation that
$$ F(z+e_n) = \lambda_{e_n} e(-m(-1)^{z_n}/4) F(z)$$
for almost all $z \in Z$.  As $m$ is odd, direct calculation then shows that $\|F(\cdot+e_n)-F\|_{L^2(Z)}$ is bounded away from zero.  But $F(\cdot+e_n)$ converges strongly to $F$ in $L^2(Z)$, giving the required contradiction.
\end{proof}

\begin{remark} One can  view $Z \rtimes_\rho K$ as a group extension of the Kronecker factor $Z \rtimes_{2\rho} 2K$ by a suitable $(Z \rtimes_{2\rho} 2K, K/2K)$-valued cocycle and obtain analogues of Proposition \ref{rho} for that cocycle; we leave the details to the interested reader.  One can also obtain higher order variants of this construction by replacing the cyclic group $\Z/4\Z$ with larger cyclic groups $\Z/2^k\Z$, or even with the $2$-adic group $\Z_2$ to create larger systems $Z \rtimes_{\rho_k} \Z/2^k\Z$ and $Z \rtimes_{\rho_\infty} \Z_2$,with the cocycles $\rho_k, \rho_\infty$ defined as in \eqref{rhog} but taking values now in $\Z/2^k\Z$ or $\Z_2$ rather than $\Z/4\Z$.  One can then show that for any $k \geq 1$, the $k^\mathrm{th}$ Host--Kra--Ziegler factor $\mathrm{Z}^k( Z \rtimes_{\rho_\infty} \Z_2 )$ is isomorphic to $Z \rtimes_{\rho_k} \Z/2^k\Z$; we leave the details of this computation to the interested reader.  Note that the previous calculations are consistent with the $k=1,2$ cases of this assertion.
\end{remark}

\begin{remark}  Essentially the same system was also studied in \cite[Appendix E]{tz-lowchar}, as an example of a system in which polynomials did not have roots of the expected degree.
\end{remark}

\subsection{Second example: a system associated to a bilinear form in odd characteristic}

Let $p$ be an odd prime, $\Gamma\coloneqq \F_p^\omega$, $Z\coloneqq\F_p^\N$, 
$\nu\coloneqq \mu_{\F_p^\N}$ its Haar measure, and $B:\Gamma\times\Gamma\to \F_p$ the standard bilinear form $B(\gamma,\gamma')\coloneqq\sum_n \gamma_n\gamma_n'$. 
We define the rotational $\Gamma$-system $\mathrm{Z}=(Z,\nu,S)$ where $S:\Gamma\to\Aut(Z,\nu)$ is the $\Gamma$-rotation $S^\gamma (z_n)\coloneqq (z_n +2\gamma_n)$.  Since $\{2 \gamma \colon \gamma\in\Gamma\}$ is dense in $Z$ (identifying $\Gamma$ with a subgroup of $Z$), the rotational $\Gamma$-system $\mathrm{Z}$ is ergodic. Let $K=\F_p$ and $\rho = (\rho_\gamma)_{\gamma \in \Gamma}$ be the $(Z, K)$-cocycle 
\begin{equation}\label{rhoh}
\rho_\gamma(z) \coloneqq \sum_{n \in \N} z_n \gamma_n + B(\gamma,\gamma)
\end{equation}
It is not difficult to verify that this is indeed a $(Z,K)$-cocycle.

We claim that $\rho$ obeys the properties stated in Proposition \ref{rho}.  We begin with ergodicity.  If this cocycle was not ergodic, then by repeating the proof of Proposition \ref{rho}(i) we could find $F \in M(Z,\T)$ such that
$$ \frac{1}{p} \rho_\gamma = F \circ V_\gamma - F \mod 1$$
$\nu$-almost everywhere for all $\gamma$.  In particular
$$ F(z + e_n) = F(z) + \frac{z_n+1}{p} \mod 1 $$
for any generator $e_n$ of $\F_p^\omega$, and this again contradicts Lusin's theorem.

To see that $\rho$ is a type $2$ cocycle, we can directly verify that \eqref{type-2} holds (with $\gamma$ replaced by $2\gamma$ on the right-hand side) with $F \equiv 0$.  Similarly, to verify the Conze--Lesigne equation, direct calculation shows that \eqref{cleq} holds with $F_z \equiv 0$ and $c_z(\gamma) \coloneqq \sum_{n \in \N} z_n \gamma_n$.

A modification of the proof of Proposition \ref{kroncalc} (relying on the assumption that $p$ is prime) reveals that $\mathrm{Z}$ is the Kronecker factor of $Z\rtimes_\rho K$; we leave the details to the interested reader.  We also can compute the Host--Kra group similarly to the previous example. Indeed, following similar computations, we find that if $(u,F)\in G$, then 
\[
F(z)= \theta + \sum_{n\in \N} z_n \sigma_n 
\]
almost surely for some $\theta\in K$ and $\sigma\in \Gamma=\F_p^\omega$. A difference to the previous example is that these $F$ correspond to the eigenfunctions of the Kronecker factor $Z$. 
Proceeding as with the previous example, we can then identify $G$ with the set of all triples $(u,\theta,\sigma)\in Z\times K\times \Gamma$ endowed  with the group law
\[
(u,\theta,\sigma)(u',\theta',\sigma')=(u+u', \theta+\theta'+\sum_{n \in \N} u_n \sigma_n, \sigma+\sigma')
\]
The remaining analysis can be carried out analogously to the previous example, and we leave it to the interested reader. 

\subsection{Third example: a system associated to a bilinear form in characteristic zero}

We now present a ``characteristic zero'' variant of the previous example, which is a standard skew-shift system.
Let $\Gamma=\Z$, $K=Z=\T$ both equipped with Lebesgue measure, and $\alpha\in \T$ be irrational.   
We equip $Z$ with the rotational $\Gamma$-system $z\mapsto z+2\alpha$, and denote the resulting system by $\mathrm{Z}$.  
For $n\in\Z$, let $\rho_n\colon \T\to \T$ be defined by $\rho_n(z)=zn+\alpha n^2$. 
Then $\rho=(\rho_n)_{n\in \Z}$ is a $(Z,K)$-cocycle. 
Now form the skew product $\mathrm{X}=\mathrm{Z}\rtimes_\rho K$. 
As $2\alpha$ is an irrational rotation, $\mathrm{Z}$ is an ergodic $\Gamma$-rotational system.
By irrationality of $\alpha$ and the definition of $\rho$, we conclude that $\mathrm{Z}$ is the Kronecker factor of $\mathrm{X}$.  
For $x=(0,0)\in \T^2$, we put $x(n)\coloneqq T_\rho^n(x)=(2\alpha n, \alpha n^2)$. 
By Weyl's equidistribution theorem (e.g., see \cite[Corollary 1.1.9]{tao2012higher}), $(x(n))_{n\in \Z}$ is asymptotically equidistributed in $\T^2$ with respect to Haar measure. 
Thus $\mathrm{X}$ is an ergodic $\Gamma$-system. 
The cocycle $\rho$ satisfies the type $2$ condition \eqref{type-2} with $F \equiv 0$, and the Conze-Lesigne equation \eqref{cleq} with $c_u(n) \coloneqq nu \mod 1$ for all $n \in \Gamma$ and $u \in Z$.  All eigenvalues of $(Z,2\alpha)$ are of the form $2\alpha m$ for $m\in \Z$. 
Therefore, we can identify the Host-Kra group $G$ of $\mathrm{X}$ with the set of all triples $(u,\theta,m)\in Z\times K\times \Gamma$ with group law 
\[
(u,\theta,m)(u',\theta',m') \coloneqq (u+u',\theta+\theta' +u, m +m')
\]
and inverse 
\[
(u,\theta,m)^{-1} \coloneqq (-u,-\theta -u, -m)
\]
We observe that $G$ acts continuously and transitively on $Z\times K$ by 
\[
(u,\theta,m)\cdot (z,k)=(z+u,k+\theta+mz). 
\]
The stabilizer of $x=(0,0)\in Z\times K$ is just 
\[
\Lambda=\{(0,0,m)\in m\in \Z\}
\]
We define 
\[
G_2 \coloneqq [G,G]=\{(0,k,0)\in k\in K\}
\]
Finally, define the translation $\phi\colon \Z\to G$ by $\phi(n) \coloneqq (2\alpha n,\alpha n^2,n)$. 
One now verifies that Theorem \ref{main-thm} holds for this example as well. 

\begin{remark}
The last two examples can be unified into a class of examples of Conze--Lesigne systems which we sketch in the following. 
Let $\Gamma$ be a countable abelian group, $U$ be a compact abelian group, and $B \colon \Gamma \times \Gamma\to U$ be a symmetric bilinear form. Let $Z=\Hom(\Gamma,U)$ (which is a compact abelian group equipped with Haar measure), and let $X$ be the set of quadratic functions $x \colon \Gamma\to U$ defined by $x(\gamma)=B(\gamma,\gamma) + z(\gamma)+c$ for some $z\in Z$ and $c\in U$. We can identify $X$ with $Z \times U$ and equip $X$ with product of Haar measures.  
We let $\Gamma$ act on $X$ by 
\[
(\gamma\cdot x)(\gamma')\coloneqq B(\gamma+\gamma',\gamma+\gamma') + z(\gamma+\gamma')+c.
\]
This extends a translational action on $Z$ defined by
$$ (\gamma \cdot z)(\gamma') \coloneqq z(\gamma) + 2 B(\gamma,\gamma')$$
using the $(Z,U)$-cocycle
$$ \rho_\gamma(z) \coloneqq B(\gamma,\gamma) + z(\gamma).$$
Under a suitable genericity hypothesis\footnote{It appears tentatively that the correct genericity hypothesis to make here is that there does not exist a finite index subgroup $\Gamma'$ of $\Gamma$ and a non-trivial character $\xi \in \hat U$ such that $\xi \circ B$ vanishes on $\Gamma' \times \Gamma'$, although we will not establish this here.}, these actions are ergodic.  The verifications of the type $2$ property \eqref{type-2} and Conze--Lesigne equation
\eqref{cleq} for $\rho$ proceed similarly to before, and one can express this system as a translational system $G/\Lambda$ with $G$ the Host--Kra group.  We leave the details to the interested reader.
\end{remark}

\begin{remark} In all of the above examples, the Host--Kra group $G$ ends up being a semi-direct product of $Z$ and $K \times \Lambda$.  However, this need not be the case in general, particularly when the cocycle $\rho$ is not of a polynomial nature.  Suppose for instance we take the Heisenberg nilsystem $G/\Lambda$ with $\Gamma \coloneqq \Z$,
$$ G \coloneqq \begin{pmatrix}
1 & \R & \R/\Z \\ 0 & 1 & \R \\ 0 & 0 & 1  \end{pmatrix}; \quad \Lambda \coloneqq \begin{pmatrix} 1 & \Z & 0 \\ 0 & 1 & \Z \\ 0 & 0 & 1 \end{pmatrix}$$
with the group action $\phi \colon \Gamma \to G$ given by
\begin{equation}\label{trans}
 \phi(n) \coloneqq \begin{pmatrix} 1 & n \alpha & \frac{n(n-1)}{2} \alpha \beta \mod 1 \\ 0 & 1 & n\beta \\ 0 & 0 & 1 \end{pmatrix}
 \end{equation}
 for some real numbers $\alpha,\beta$ with $1,\alpha,\beta$ linearly independent over the rationals.  The Kronecker factor $Z$ can be identified with the two-torus $(\R/\Z)^2$ with translation map $S^n \colon (x,y) \mapsto (x+n\alpha, y+n\beta)$, and by following the construction in Section \ref{nil-cl} with the section $s \colon Z \to G/\Lambda$ defined by
$$ s(x,y) \coloneqq \begin{pmatrix} 1 & \{x\} & 0 \\ 0 & 1 & \{y\} \\ 0 & 0 & 1 \end{pmatrix}$$
with $x \mapsto \{x\}$ the fractional part map from $\R/\Z$ to $[0,1)$, we can calculate the cocycle $(Z,K)$-cocycle $\rho$ (with $K = \R/\Z$) to be
$$\rho_n(x,y) = \frac{n(n-1)}{2} \alpha \beta + n \alpha \{y\} - (x + n\alpha) (\{y\}+n\beta - \{y+n\beta\} ) \mod 1.
$$
Here the Host--Kra group $G$ is not the semi-direct product of $Z$ and $K \times \Lambda$; instead we have a non-split short exact sequence
$$ 0 \to H \to G \to Z \to 0$$
with 
$$ H \coloneqq \begin{pmatrix} 1 & \Z & \R/\Z \\ 0 & 1 & \Z \\ 0 & 0 & 1 \end{pmatrix} \equiv K \times \Lambda.$$
\end{remark}

\begin{remark}[Rudolph's example]\label{rudolph-ex}\footnote{This reformulation of Rudolph's example was communicated to us by Yonatan Gutman.}  Let $\alpha,\beta$ be as in the previous example.  One can take an inverse limit of the $\Z$-nilsystems
$$  \begin{pmatrix}
1 & \R & \R/2^N \Z \\ 0 & 1 & \R \\ 0 & 0 & 1  \end{pmatrix} / \begin{pmatrix} 1 & 2^N \Z & 0 \\ 0 & 1 & \Z \\ 0 & 0 & 1 \end{pmatrix}$$
as $N \to \infty$, using the translation action \eqref{trans} for each $N$, to obtain a Conze--Lesigne system that was shown by Rudolph \cite{rudolph} to not be expressible as a nilpotent translational $\Z$-system of nilpotency class two.  It can be expressed as an abelian extension $Z \rtimes_\rho K$, where the Kronecker factor $Z$ is given by $Z \coloneqq S_2 \times \R/\Z$, with $S_2$ is the $2$-adic solenoid group formed as the inverse limit of the $\R/2^N\Z$, and where the translation action $n \mapsto (n \alpha, n \beta)$, $K \coloneqq S_2$ is another copy of the solenoid group, and the $(Z,K)$-cocycle $\rho$ is given by
$$ \rho_n(x,y) \coloneqq \frac{n(n-1)}{2} \alpha \beta + n \alpha \{y\} - (x + n\alpha) (\{y\}+n\beta - \{y+n\beta\} ) $$
(where we embed $\R$ into $S_2$ in the obvious fashion) for $n \in \Z$, $x \in S_2$, $y \in \R/\Z$, noting that the product in the last term is well-defined since $\{y\}+n\beta - \{y+n\beta\}$ is an integer.  This cocycle is ergodic and of type $2$ but does not obey the Conze--Lesigne equation, mainly because there are too few continuous homomorphisms\footnote{In particular, there are no non-trivial continuous homomorphisms from $\T$ to $S_2$.} from $Z$ to  $K$; this does not contradict Theorem \ref{conze-eq} because $K$ is not a Lie group.  On the other hand, the system $Z \rtimes_\rho K$ can be expressed as a double coset system
$$ \begin{pmatrix} 1 & 0 & 0 \\ 0 & 1 & \{0\} \times \Z_2 \\ 0 & 0 & 1 \end{pmatrix} \backslash \begin{pmatrix}
1 & \R \times \Z_2 & \R \times \Z_2 \\ 0 & 1 & \R \times \Z_2 \\ 0 & 0 & 1  \end{pmatrix} / \begin{pmatrix} 1 & \Delta(\Z) & \Delta(\Z) \\ 0 & 1 & \Delta(\Z) \\ 0 & 0 & 1 \end{pmatrix}$$
where the $2$-adic group $\Z_2$ is the inverse limit of $\Z/2^N\Z$ and $\Delta$ is the diagonal embedding of $\Z$ into $\R \times \Z_2$; see \cite{shalom2}.
\end{remark}

\appendix

\section{Concrete and abstract measure theory}\label{conc-app}

In this appendix we review the notational conventions we will use for various types of probability spaces, and measure-preserving actions on such spaces.  It will be convenient to use some of the category theoretic formalism from \cite{jt-foundational}, although we will not make heavy use of category-theoretic tools in this paper.

\subsection{Forgetful functors}\label{forget}

We begin with a general convention concerning ``casting functors'' from \cite{jt-foundational}, although in this paper we will refer to these functors as ``forgetful functors'' instead.

We will deem a number of functors
\begin{center}
\begin{tikzcd}
    \Cat \arrow[r,hook] & \Cat_0
\end{tikzcd}
\end{center}
from one category $\Cat$ to another $\Cat_0$ to be ``forgetful functors'', which intuitively would take an $\Cat$-object $X$ or a $\Cat$-morphism $f \colon X \to Y$ and ``forget'' some of its structure to return a $\Cat_0$-object $X_{\Cat_0}$ or a $\Cat_0$-morphism $f_{\Cat_0} \colon X_{\Cat_0} \to Y_{\Cat_0}$.  We always consider the identity functor to be forgetful, and the composition of two forgetful functors to be forgetful\footnote{This convention will be unambiguous because all of our forgetful functors will commute with each other.}; for instance if we have two forgetful functors
\begin{center}
\begin{tikzcd}
    \Cat_1 \arrow[r,hook] & \Cat_2  \arrow[r,hook] & \Cat_3
\end{tikzcd}
\end{center}
then we have $X_{\Cat_3} = (X_{\Cat_2})_{\Cat_3}$ for any $\Cat_1$-object $X$.

Given a pair of forgetful functors
\begin{center}
\begin{tikzcd}
    \Cat_1 \arrow[r,hook] & \Cat_0 & \Cat_2 \arrow[l,hook'] 
\end{tikzcd}
\end{center}
we say that a $\Cat_1$-object $X_1$ and a $\Cat_2$-object $X_2$ are \emph{$\Cat_0$-isomorphic} if there is a $\Cat_0$-isomorphism between $(X_1)_{\Cat_0}$ and $(X_2)_{\Cat_0}$.  Similarly, a $\Cat_1$-morphism $f_1 \colon X_1 \to Y_1$ and a $\Cat_2$-morphism $f_2 \colon X_2 \to Y_2$ are \emph{$\Cat_0$-equivalent} if $X_1,X_2$ are $\Cat_0$-isomorphic, $Y_1,Y_2$ are $\Cat_0$-isomorphic, and the $\Cat_0$-morphisms $(f_1)_{\Cat_0} \colon (X_1)_{\Cat_0} \to (Y_1)_{\Cat_0}$, $(f_2)_{\Cat_0} \colon (X_2)_{\Cat_0} \to (Y_2)_{\Cat_0}$ agree after composing with these $\Cat_0$-isomorphisms.  If a $\Cat_1$-object $X_1$ is $\Cat_0$-isomorphic to a $\Cat_0$-object $X'$, we call $X_1$ a \emph{$\Cat_1$-model} of $X'$; similarly, if a $\Cat_1$-morphism $f_1 \colon X_1 \to Y_1$ is $\Cat_0$-equivalent to a $\Cat_0$-morphism $f' \colon X' \to Y'$, we call $f_1$ a \emph{$\Cat_1$-representation} of $f'$.

\subsection{Probability spaces}

In this paper we will work with three categories $\ConcProb, \ProbAlg, \CHProb$ of probability spaces.

\begin{definition}[Categories of probability spaces]\label{prob-def}\cite{jt-foundational}
\begin{itemize}
    \item[(i)]  A \emph{concrete probability space} (or \emph{$\ConcProb$-space}) is a triple $(X,\X,\mu)$, where $X$ is a set, $\X$ is a $\sigma$-algebra of subsets of $X$, and $\mu \colon \X \to [0,1]$ is a countably additive probability measure.  A \emph{concrete probability-preserving map} (or \emph{$\ConcProb$-morphism}) $f \colon (X,\X,\mu) \to (Y,\Y,\nu)$ between two $\ConcProb$-spaces is a measurable map $f \colon X \to Y$ such that $\mu(f^{-1}(F)) = \nu(F)$ for all $F \in \Y$ (that is to say, the pushforward $f_* \mu$ of $\mu$ by $f$ is equal to $\nu$). 
    \item[(ii)]  A \emph{probability algebra} (or $\OpProbAlg$-space) is an object of the form $(\X,\mu)$, where $\X = (\X, \vee, \wedge, 0, 1, \bar{\cdot})$ is a $\sigma$-complete Boolean algebra, and $\mu \colon \X \to [0,1]$ is a countably additive probability measure on $\X$ such that $\mu(E)=0$ if and only if $E=0$.   An \emph{abstract probability-preserving map} (or \emph{$\OpProbAlg$-morphism}) $f \colon (\X,\mu) \to (\Y,\nu)$ between two $\OpProbAlg$-spaces is a Boolean homomorphism\footnote{Notice the opposite direction of the arrows here. We implicitly work with an opposite category here to keep certain functors covariant.} $f \colon \Y\to \X$ that preserves countable joins\footnote{Actually, the preservation of countable joins is automatic for Boolean homomorphisms between probability algebras, and such algebras are in fact complete Boolean algebras as opposed to merely being $\sigma$-complete, although we will not need these (easily established) facts here.} (thus $f(\bigvee_{n=1}^\infty F_n) = \bigvee_{n=1}^\infty f(F_n)$ for all $F_n \in \X$ such that $\mu(f(F)) = \nu(F)$ for all $F \in \X$).
    \item[(iii)]  A \emph{compact probability space} (or \emph{$\CHProb$-space}) is a quadruple $(X,{\mathcal F}, \X,\mu)$, where $(X,{\mathcal F})$ is a compact Hausdorff topological space, $\X$ is the Baire\footnote{See \cite{jt19, jt-foundational} for a discussion as to why the Baire $\sigma$-algebra is a more natural choice than the Borel $\sigma$-algebra for compact Hausdorff spaces that are not necessarily metrizable, and similarly for why ``compact $G_\delta$ inner regular in the Baire algebra'' is the natural definition of a Radon measure in this setting.} $\sigma$-algebra (i.e., the topology generated by the continuous functions from $X$ to $\R$), and $\mu$ is a countably additive probability measure which is Radon in the sense of \cite[Definition 4.1]{jt-foundational}, i.e., $\mu$ is compact $G_\delta$ inner regular in the Baire algebra.  A \emph{continuous probability-preserving map} (or \emph{$\CHProb$-morphism}) $f \colon (X,{\mathcal F},\X,\mu) \to (Y,{\mathcal G}, \Y,\nu)$ between $\CHProb$-spaces is a continuous map which is also a $\ConcProb$-morphism.
\end{itemize}
\end{definition}

It is easy to verify that $\ConcProb, \OpProbAlg, \CHProb$ are indeed categories.  Inside these categories we isolate some ``countable'' objects:

\begin{itemize}
    \item[(i)]  A concrete probability space $(X,\X,\mu)$ is a \emph{Lebesgue space} (or \emph{Lebesgue} for short) if the measurable space $(X,\X)$ is a standard Borel space, that is to say one can endow $X$ with the structure of a Polish space such that $\X$ is the Borel $\sigma$-algebra.  
    \item[(ii)] A probability algebra $(\X,\mu)$ is \emph{separable} if the $\sigma$-complete Boolean algebra $\X$ is countably generated.
    \item[(iii)] A compact probability space $(X,{\mathcal F}, \X,\mu)$ is \emph{metrizable} if the topological space $(X,{\mathcal F})$ is metrizable (or equivalently by the Urysohn metrization theorem, second countable).  
\end{itemize}

There are obvious \emph{forgetful functors}
\begin{center}
\begin{tikzcd}
    \CHProb \arrow[r,hook] & \ConcProb \arrow[r,hook] & \OpProbAlg,
\end{tikzcd}
\end{center}
between these categories, in which a $\CHProb$-space $(X,{\mathcal F},\X,\mu)$ is converted to a $\ConcProb$-space $(X,{\mathcal F},\X,\mu)_\ConcProb \coloneqq (X,\X,\mu)$ by forgetting the topology ${\mathcal F}$, and a $\ConcProb$-space $(X,\X,\mu)$ is converted to a probability algebra $(X,\X,\mu)_\OpProbAlg \coloneqq (\X_\mu, \mu)$ by forming the \emph{probability algebra}
$$ \X_\mu \coloneqq \{ [E] \colon E \in \X \}$$
where for each $E \in \X$, the equivalence class $[E]$ is defined as the collection of sets equal modulo null sets to $E$, thus
$$ [E] \coloneqq \{ F \in \X \colon \mu(E \Delta F) = 0\},$$
and by abuse of notation we define $\mu \colon \X_\mu \to [0,1]$ by requiring $\mu([E]) \coloneqq \mu(E)$ for all $E \in \X$.  Morphisms are then also transformed in the obvious fashion (although the direction of the arrows are "flipped" when moving from $\ConcProb$ to $\OpProbAlg$).  It is not difficult to verify that these are indeed functors, and we will adopt the forgetful functor conventions from Section \ref{forget}.  We also describe some of these conventions in plainer English:

\begin{itemize}
    \item If $X$ is a compact or concrete probability space, we refer to the probability algebra  $X_\OpProbAlg$ as the \emph{abstraction}\footnote{More precisely this should be ``abstraction modulo null sets'', as we are both abstracting away the space $X$ and quotienting out by the null ideal. Similarly for other uses of the term ``abstract'' in this paper.} of $X$.  Similarly, if $f \colon X \to Y$ is a continuous or concrete probability-preserving map, we refer to the abstract probability-preserving map $f_\OpProbAlg \colon Y_\OpProbAlg \to X_\OpProbAlg$ as the \emph{abstraction} of $f$.
    \item A $\ConcProb$-model $\tilde X$ (resp. $\CHProb$-model $\hat X$) of a probability algebra $X$ will be called a \emph{concrete model} (resp. \emph{topological model}) of $X$.  Similarly, a $\ConcProb$-representation $\tilde f \colon \tilde X \to \tilde Y$ (resp. $\CHProb$-representation  $\hat f \colon \hat X \to \hat Y$) of an abstract probability-preserving map $f \colon Y \to X$ will be called a \emph{concrete representation} (resp. \emph{continuous representation}) of $f$.
\end{itemize}

Observe that if a compact probability space $(X,{\mathcal F}, \X,\mu)$ is metrizable, then the Baire $\sigma$-algebra coincides with the Borel $\sigma$-algebra and so the associated concrete probability space $(X,{\mathcal F}, \X,\mu)_\ConcProb$ is Lebesgue.  Similarly, if a concrete probability space $(X,\X,\mu)$ is Lebesgue, then the associated probability algebra $(X,\X,\mu)_\OpProbAlg$ is separable.  Thus the notions of ``countability'' for the three categories $\ConcProb$, $\OpProbAlg$, $\CHProb$ are all compatible with each other.  On the other hand, the converse implications are false; it is entirely possible for a separable probability algebra to be modeled by a concrete probability space that is not Lebesgue, or a compact probability space that is not metrizable.

If two concrete measure-preserving maps $f, g \colon X \to Y$ agree almost everywhere, then they are abstractly equal: $f_\OpProbAlg = g_{\OpProbAlg}$.  However, if the target space $Y$ is not Lebesgue or Polish, the converse statement can fail; see \cite[Examples 5.1, 5.2]{jt19}.  Nevertheless the reader may wish to think of ``agreement almost everywhere'' as a heuristic first approximation of the concept of ``abstract equality''.

It is natural to ask to what extent the above forgetful functors can be inverted.  In this regard we have the following results:

\begin{proposition}[Reversing the forgetful functors for probability spaces]\label{reverse}\ 
\begin{itemize}
    \item[(i)] (Existence of concrete representations) \cite[Proposition 3.2]{jt19} If $(X,\X,\mu)$, $(Y,\Y,\nu)$ are concrete probability spaces  with $(Y,\Y,\nu)$ Lebesgue, then every abstract probability-preserving map $f \colon (X,\X,\mu)_\OpProbAlg \to (Y,\Y,\nu)_\OpProbAlg$ has a concrete representation $\tilde f \colon (Y,\Y,\nu)  \to (X,\X,\mu)$, which is a concrete probability-preserving map that is unique up to almost everywhere equivalence.  Related to this, two concrete measurable maps from $X$ to a Polish space $Y$ agree abstractly if and only if they agree almost everywhere.
    \item[(ii)] (Cantor model) \cite[Theorem 2.15]{glasner2015ergodic}  If $\pi \colon (\X,\mu) \to (\Y,\nu)$ is an abstract probability-preserving map between separable probability algebras, then there exists a continuous representation $\hat \pi \colon (\hat Y,{\mathcal G},\hat \Y,\hat \nu)\to (\hat X,{\mathcal F},\hat \X,\hat \mu) $ of $\pi$ between compact metrizable probability spaces (in fact Cantor spaces).
    \item[(iii)] (Canonical model) \cite[Theorem 7.2]{jt-foundational} There exists a \emph{canonical model functor} (or \emph{Stone functor})
\begin{center}
\begin{tikzcd}
    \OpProbAlg \arrow[r,"\Stone", hook, two heads] & \CHProb, 
\end{tikzcd}
\end{center}
that takes a probability algebra $(\X,\mu)$ and constructs a topological model $\Stone(\X,\mu)$, and similarly takes any opposite abstract probability-preserving map $f \colon (\X,\mu) \to (\Y,\nu)$ and constructs (in a completely functorial and natural fashion) a continuous representation  $\Stone(f) \colon \Stone(\X,\mu) \to \Stone(\Y,\nu)$.
\end{itemize}
\end{proposition}

We remark that the compact probability space $\Stone(\X,\mu)$ in Proposition \ref{reverse}(iii) is constructed using either Gelfand duality or Stone duality and is not metrizable in general, even when $(\X,\mu)$ is separable.

Given a concrete probability space $(X,\X,\mu)$, we let $L^0(X,\X,\mu)$ denote the space of measurable functions from $X$ to $\C$, quotiented out by almost everywhere equivalence, and for $1 \leq p \leq \infty$ we let $L^p(X,\X,\mu)$ denote the subspace of $L^0(X,\X,\mu)$ consisting of those (equivalence classes of) measurable functions which are $p^{\mathrm{th}}$ power integrable (or essentially bounded, in the $p=\infty$ case).  For a probability algebra $(\X,\mu)$ we can similarly define $L^0(\X,\mu)$ and $L^p(\X,\mu)$ by passing to a concrete or topological model (for instance by using the canonical model functor $\Stone$) and using the previous construction; note that up to isomorphism, the precise choice of model used is irrelevant.  Note that $L^2(\X,\mu)$ is a Hilbert space and $L^\infty(\X,\mu)$ is a tracial commutative von Neumann algebra (and hence also a $C^*$-algebra), using the integral against $\mu$ as the trace.

\subsection{Dynamics}\label{dynamic}

Let $\Gamma$ be an arbitrary group (not necessarily countable or abelian); for this discussion we treat $\Gamma$ as a discrete group, ignoring any topological structure.  Let $\Cat$ be one of the three categories $\ConcProb$, $\OpProbAlg$, $\CHProb$.  We define a \emph{$\Cat_\Gamma$-system} to be a pair $\mathrm{X} = (X,T)$, where $X$ is a $\Cat$-space and $T \colon \Gamma \to \Aut_\Cat(X)$ is a group homomorphism of $\Gamma$ to the automorphism group $\Aut_\Cat(X)$, that is to say the group of $\Cat$-isomorphisms from $X$ to itself.  A \emph{$\Cat_\Gamma$-morphism} $\pi \colon (X,T) \to (Y,S)$ between two $\Cat_\Gamma$-systems $(X,T)$, $(Y,S)$ is a $\Cat$-morphism $\pi \colon X \to Y$ with the property that one has the identity $S^\gamma \circ \pi = \pi \circ T^\gamma$ of $\Cat$-morphisms for all $\gamma \in \Gamma$.  Properties defined for $\Cat$-spaces are then also applicable to $\Cat_\Gamma$-systems in the obvious fashion; for instance, a $\ConcProb_\Gamma$-system $(X,T)$ is Lebesgue if the underlying $\ConcProb$-space $X$ is Lebesgue.  We also adopt the following terminology:

\begin{itemize}
    \item $\ConcProb_\Gamma$-systems and $\ConcProb_\Gamma$-morphisms will be called \emph{concrete $\Gamma$-systems} and \emph{concrete factor maps} respectively.
    \item $\OpProbAlgG$-systems and $\OpProbAlgG$-morphisms will be called \emph{abstract $\Gamma$-systems} and \emph{abstract factor maps} respectively.  Any two $\Gamma$-systems will be called \emph{abstractly isomorphic} if they are $\OpProbAlgG$-isomorphic.
    \item $\CHProb_\Gamma$-systems and $\CHProb_\Gamma$-morphisms will be called \emph{compact $\Gamma$-systems} and \emph{continuous factor maps} respectively.
\end{itemize}
    
The diagram of forgetful functors from the previous subsection can now be enlarged to a commuting diagram
\begin{center}
\begin{tikzcd}
    \CHProb_\Gamma \arrow[r,hook] \arrow[d,hook] & \ConcProb_\Gamma \arrow[r,hook] \arrow[d,hook] & \OpProbAlgG \arrow[d,hook] \\
    \CHProb \arrow[r,hook] & \ConcProb \arrow[r,hook] & \OpProbAlg
\end{tikzcd}
\end{center}
of forgetful functors in the obvious fashion.  We adapt concepts such as topological models, concrete representations, etc. to this dynamical setting; for instance, a $\CHProb_\Gamma$-model $\hat{\mathrm{X}}$ of a $\OpProbAlgG$-system $\mathrm{X}$ will be referred to as a \emph{topological model} of $\mathrm{X}$.

For us, one important source of an abstract $\Gamma$-system arises by starting with a concrete probability space $(X,\X,\mu)$ and equipping it with a \emph{near-action}\footnote{Here we follow the notation of Zimmer \cite{zimmer}.} of $\Gamma$, by which we mean a family of concrete measure-preserving maps $T^\gamma \colon X \to X$ for each $\gamma \in \Gamma$ such that $T^1(x)=x$ for $\mu$-almost all $x \in X$, and $T^{\gamma_1} T^{\gamma_2}(x) = T^{\gamma_1 \gamma_2}(x)$ for all $\gamma_1,\gamma_2 \in \Gamma$ and $\mu$-almost all $x \in X$ (with the obvious changes if the group $\Gamma$ is written additively instead of multiplicatively).  This is not quite a concrete $\Gamma$-system because of the possibility that the identities $T^1(x)=x$, $T^{\gamma_1} T^{\gamma_2}(x) = T^{\gamma_1 \gamma_2}(x)$ fail on a null set.  However, by passing to the abstract setting we see that $(T^1)_{\OpProbAlg}$ is the identity and $(T^{\gamma_1})_{\OpProbAlg}(T^{\gamma_2})_{\OpProbAlg} = T^{\gamma_1 \gamma_2}_{\OpProbAlg}$ for all $\gamma_1,\gamma_2 \in \Gamma$, so the near-action induces an abstract $\Gamma$-system
$((X,\X,\mu)_\OpProbAlg, (T^\gamma_\OpProbAlg)_{\gamma \in \Gamma})$.
    
If one has an abstract factor map $\pi \colon \mathrm{X}_{\OpProbAlgG} \to \mathrm{Y}_{\OpProbAlgG}$ between two (concrete, abstract, or compact) $\Gamma$-systems $\mathrm{X}, \mathrm{Y}$, we write $\mathrm{Y} \leq \mathrm{X}$; this is a partial order up to abstract isomorphism.  This map generates a \emph{factor algebra} $\{ \pi^* E: E \in \Y \} \subset \X$, where $\X$, $\Y$ are the $\sigma$-complete boolean algebras associated to $\mathrm{X}$, $\mathrm{Y}$ respectively.
 A factor $\mathrm{Y}$ of an abstract $\Gamma$-system $\mathrm{X}$ is said to be the \emph{inverse limit} of a collection $(\mathrm{Y}_\alpha)_{\alpha \in A}$ of factors indexed by a directed set $A$ (with factor maps $\pi_{\alpha \beta} \colon \mathrm{Y}_\beta \to \mathrm{Y}_\alpha$ whenever $\alpha\leq \beta$ that all commute with each other and with the factor maps $\pi_\alpha: \mathrm{Y} \to \mathrm{Y}_\alpha$ in the obvious fashion) if the factor algebra of $\mathrm{Y}$ is generated by the union of the factor algebras of the $\mathrm{Y}_\alpha$.

Again, we have some results concerning the extent to which the forgetful functors can be inverted:

\begin{proposition}[Reversing the forgetful functors for systems]\label{reverse-system}  Let $\Gamma$ be a group.
\begin{itemize}
\item[(i)]  (Concrete representation) \cite[Theorem 2.15(ii)]{glasner2015ergodic}  If $\Gamma$ is countable, $\mathrm{X}, \mathrm{X}'$ are concrete Lebesgue $\Gamma$-systems, and $\pi \colon \mathrm{X}_{\OpProbAlgG} \to \mathrm{X}'_{\OpProbAlgG}$ is an isomorphism of abstract $\Gamma$-systems, then there exist full measure concrete sub-systems $\mathrm{X}_0, \mathrm{X}'_0$ of $\mathrm{X}, \mathrm{X}'$ (formed by deleting $\Gamma$-invariant Borel null sets from both systems) and a concrete representation $\tilde \pi \colon \mathrm{X}_0 \to \mathrm{X}'_0$ of $\pi$.
\item[(ii)]  (Cantor representation) \cite[Theorem 2.15(i)]{glasner2015ergodic}  If $\Gamma$ is countable, and $\pi \colon \mathrm{X} \to \mathrm{Y}$ is an abstract factor map between abstract separable $\Gamma$-systems, then there exists a continuous representation $\hat \pi \colon \hat{\mathrm{X}} \to \hat{\mathrm{Y}}$ of $\pi$ between compact metrizable $\Gamma$-systems (in fact Cantor systems).
\item[(iii)] (Canonical model) \cite[Theorem 7.2]{jt-foundational} The canonical model functor from Proposition \ref{reverse}(iii) induces a commuting square of functors
\begin{center}
\begin{tikzcd}
    \OpProbAlgG \arrow[r,"\Stone", hook, two heads] \arrow[d,hook] & \CHProb_\Gamma \arrow[d,hook]  \\
    \OpProbAlg \arrow[r,"\Stone", hook, two heads] & \CHProb 
\end{tikzcd}
\end{center}
in the obvious fashion, such that if $\pi \colon \mathrm{X} \to \mathrm{Y}$ is an abstract factor map between abstract $\Gamma$-systems, then the continuous factor map  $\Stone(\pi) \colon \Stone(\mathrm{X}) \to \Stone(\mathrm{Y})$ is a continuous representation of $\pi$ (and $\Stone(\mathrm{X})$, $\Stone(\mathrm{Y})$ are topological models of $\mathrm{X}, \mathrm{Y}$ respectively).
\item[(iv)] (Concrete representation, II) \cite[Proposition 3.1]{zimmer}, \cite[Lemma 3.2]{ramsay} If $\Gamma$ is countable, $\mathrm{X}= (X,T)$ is an abstract $\Gamma$-system, and $\tilde X$ is a concrete Lebesgue model for $X$, then there exists a concrete model $\tilde{\mathrm{X}} = (\tilde X, \tilde T)$ of $\mathrm{X}$.
\end{itemize}
\end{proposition}

\subsection{Koopman models}\label{koopman-sec}

We now construct a topological model $\hat X$ that one can associate to any $\OpProbAlg$-space $X$ that has an action of a locally compact group $G$.  This model is constructed via the Koopman action of $G$ and so we refer to this as the \emph{Koopman model} of $X$; this generalizes the canonical model $\Stone(X)$ discussed earlier, which corresponds to the case when the group $G$ is trivial.  Our treatment is inspired by that in \cite[\S 19.3.1]{hk-book}, \cite{hk-errata}.  By taking advantage of the general category theoretic dualities in \cite{jt-foundational}, we can avoid the need to impose any ``countability'' or ``separability'' hypotheses on our spaces and groups.

\begin{theorem}[Koopman model]\label{koop}  Let $G$ be a group (not necessarily countable, discrete, or abelian), and let $\mathrm{X} = (X,T)$ be an abstract $G$-system.  Assume the following axioms:
\begin{itemize}
    \item[(i)] The $G$-action is \emph{abstractly faithful} in the sense that the \emph{Koopman representation} $g \mapsto U_g$, which assigns to each $g \in G$, the unitary \emph{Koopman operator} $U_g \colon L^2(X) \to L^2(X)$ defined by
$$ U_g(f) \coloneqq f \circ T_g^{-1},$$
is injective.    
    \item[(ii)] By (i),  we can identify $G$ with a subgroup of the unitary group of $L^2(X)$, Endowing the latter with the strong operator topology, we assume that $G$ is locally compact.  
\end{itemize}
Then there exists a topological model $\hat{\mathrm{X}} = (\hat X, \hat T) = (\hat X, {\mathcal F}, \hat \X, \hat \mu, \hat T)$ of $\mathrm{X} = (X,T)$ (which we call the \emph{Koopman model} of $\mathrm{X}$) with the following properties:
\begin{itemize}
    \item[(a)]  All non-empty open sets in $\hat X$ have positive measure.
    \item[(b)]  The action $\hat T \colon G \times \hat X \to \hat X$ is jointly continuous in $G$ and $\hat X$ (as opposed to merely being continuous in $\hat X$ for each individual group element $g \in G$).
    \item[(c)]  If $f \in L^\infty(X)$ is \emph{$G$-continuous} in the sense that the map $g \mapsto U_g(f)$ is a continuous map from $G$ to $L^\infty(X)$, then $f$ has a continuous representative $\hat f \in C(\hat X)$ in $\hat X$ (which is unique by property (a)).
\end{itemize}
Furthermore, the model $\hat{\mathrm{X}}$ is unique up to isomorphism of compact $G$-systems.
\end{theorem}

\begin{proof} 
We first establish uniqueness of the Koopman model $\hat{\mathrm{X}}$.  Being a topological model, we can identify $L^\infty(\hat X)$ with $L^\infty(X)$ as a tracial commutative $C^*$ algebra.  There is an obvious tracial $C^*$-algebra homomorphism from $C(\hat X)$ to $L^\infty(\hat X) \equiv L^\infty(X)$, which is injective from property (a).  From property (c), the image of this homomorphism contains all the $G$-continuous functions; conversely, from property (b), every element of this image is $G$-continuous.  Thus as a tracial commutative $C^*$-algebra, $C(\hat X)$ (viewed as a subalgebra of $L^\infty(X)$) is uniquely determined by the abstract $G$-system $\mathrm{X}$.  The uniqueness of the model up to isomorphism then follows from the Gelfand--Riesz duality (i.e., Gelfand duality combined with the Riesz representation theorem) between $\CHProb$-spaces and tracial $C^*$-algebras; see \cite[Theorem 5.11]{jt-foundational}.

We now reverse these steps to establish existence of the Koopman model.  Let ${\mathcal A}$ denote the space of $G$-continuous functions in $L^\infty(X)$.  This is clearly a tracial commutative $C^*$-algebra.  We claim that the closed unit ball of this algebra is dense in the closed unit ball of $L^\infty(X)$ in the $L^2(X)$ topology.  To see this, fix a left-invariant Haar measure $dg$ on $G$, let $f$ be in the closed unit ball of $L^\infty(X)$, and consider the convolution
$$ \phi * f \coloneqq \int_G \phi(g) U_g(f)\ dg$$
of $f$ with a continuous compactly supported function $\phi \in C_c(G)$.  As $G$ is given the strong operator topology, it is easy to see that this integral is well-defined and is $G$-continuous; also, by choosing $\phi$ to be a suitable approximation to the identity (non-negative, supported on a small neighborhood of the identity, and of total mass one) and again using the fact that the topology of $G$ is given by the strong operator topology, one can ensure that $\phi*f$ lies in the closed unit ball of ${\mathcal A}$ and is arbitrarily close to $f$ in $L^2(X)$; see \cite[\S 18.3.1, Lemma 7]{host2005nonconventional}.  This establishes density.

By Gelfand--Riesz duality \cite[Theorem 5.11]{jt-foundational}, we can now construct a $\CHProb$-space $\hat X$ such that $C(\hat X)$ is isomorphic as a tracial commutative $C^*$-algebra to ${\mathcal A}$.  Identifying these two algebras, we see that the $L^2(\hat X)$ norm on $C(\hat X)$ agrees with the $L^2(X)$ norm on ${\mathcal A}$.  In particular, every non-zero element of $C(\hat X)$ has positive $L^2(\hat X)$ norm (i.e., the trace is faithful), which gives (a) by Urysohn's lemma.  As $\hat X$ is equipped with a Radon measure, $C(\hat X)$ is dense in $L^2(\hat X)$, hence on taking $L^2$ closures of unit balls we obtain an identification of $L^\infty(\hat X)$ with $L^\infty(X)$, which one can easily verify to be an isomorphism of tracial commutative von Neumann algebras.  From this and the duality of categories between tracial commutative von Neumann algebras and probability algebras  (see \cite[Theorem 7.1]{jt-foundational}) we see that $\hat X$ is a topological model of $X$.

The claim (c) is clear from construction, so it remains to establish the claim (b).  By definition of a $\CHProb_G$-system, the action $\hat T_g \colon \hat X \to \hat X$ associated to any group element $g \in G$ is an element of the space $C(\hat X, \hat X)$ of continuous maps from $\hat X$ to itself.  We endow this space with the compact-open topology.  To prove joint continuity, it then suffices to show that the map $g \mapsto \hat T_g$ is continuous from $G$ to $C(\hat X, \hat X)$.  By the homomorphism property of the group action, it suffices to show that for any net $g_\alpha$ converging to the identity in $G$, the maps $\hat T_{g_\alpha} \colon \hat X \to \hat X$ converge to the identity in the compact-open topology.  From the identification of $C(\hat X)$ with ${\mathcal A}$ we see that for any $f \in C(\hat X)$, $f \circ \hat T_{g_\alpha}$ converges uniformly to $f$, and the claim now follows from Urysohn's lemma.
\end{proof}

\begin{remark}\label{bookfix} Even when the original $\OpProbAlg$-space $X$ is separable, the Koopman model $\hat X$ need not be metrizable if the action of the group $G$ is insufficiently ``transitive''.  For instance if $G$ is the trivial group then the Koopman model $\hat X$ coincides with the canonical model $\Stone(X)$ from Proposition \ref{reverse}(iii) (basically because all elements of $L^\infty(X)$ are $G$-continuous in this case), which as previously remarked is almost never metrizable in practice.  However, if $X$ is separable and $G$ is ``weakly transitive'' in the sense that the convolution operators $f \mapsto \phi * f$ used in the above proof map $L^2(X)$ to $L^\infty(X)$ for any $\phi \in C_c(G)$, then it is not difficult to show that the $C^*$-algebra ${\mathcal A}$ is separable, and hence the Koopman model $\hat X$ will be metrizable.  This weak transitivity property is not true for arbitrary groups $G$, but can be verified for the specific Host--Kra groups arising for instance in the proof of Theorem \ref{hk-thm}; see the erratum to \cite[Chapter 19]{hk-book} at \cite{hk-errata} for more details.
\end{remark}

For us, the main application of Koopman models is to enable one to identify abstract systems as translational systems.  We formalize this using the following lemma:

\begin{lemma}[Criterion for being isomorphic to a translational system]\label{translation}  Let $G, \mathrm{X}$ obey the axioms of Theorem \ref{koop}, and let $\hat{\mathrm{X}} = (\hat X, \hat T)$ be the Koopman model.  Let $\hat x_0$ be a point in $\hat X$. Assume the following additional axioms:
\begin{itemize}
    \item[(iii)] $G$ is unimodular.
    \item[(iv)] The action of $G$ on $\hat X$ is transitive.  That is to say, for any $\hat x_1,\hat x_2 \in \hat X$, there exists $g \in G$ such that $\hat T^g \hat x_1 = \hat x_2$.
    \item[(v)]  The stabilizer group $\Lambda \coloneqq \{ g \in G: \hat T^g \hat x_0 = \hat x_0 \}$ is a lattice in $G$.
\end{itemize}
Then the Koopman model $\hat{\mathrm{X}}$ is isomorphic as a compact $G$-system to the translational $G$-system $G/\Lambda$ (with the obvious action of $G$).  In particular, the abstract $G$-system $\mathrm{X}$ is abstractly isomorphic to $G/\Lambda$.
\end{lemma}

\begin{proof}  By axiom (iv), we can form a bijection between $G/\Lambda$ and $\hat X$ by identifying any coset $g\Lambda$ with $\hat T^g x_0$.  From Theorem \ref{koop}(b) and axiom (v) this bijection is continuous; since $G/\Lambda$ and $\hat X$ are both compact Hausdorff spaces, this bijection is thus a homeomorphism, and so we may identify $G/\Lambda$ and $\hat X$ as compact Hausdorff spaces at least.  By construction, the action of $G$ on both these spaces agree; by the uniqueness of Haar probability measure on $G/\Lambda$ (which is well-defined by axioms (iii), (v)) we conclude that the measure $\hat \mu$ on $\hat X$ agrees with the Haar probability measure on $G/\Lambda$.  The claim follows.
\end{proof}

\section{Measurable selection lemma}\label{mes-select-app}

In this appendix we give a full proof of Proposition \ref{mes-select}.  The arguments here can also be used to give a more detailed proof of \cite[Lemma C.4]{btz}; we leave this modification to the interested reader.

Let $\Gamma, \mathrm{Y}, U, h_u$ be as in the proposition.  For each $u \in U$, we introduce the set
$$ \Omega_u \coloneqq \{ F \in M(\mathrm{Y}, \T): h_u - dF \in \hat \Gamma \},$$
then by hypothesis $\Omega_u$ is non-empty for each $u$.  Observe that each $\Omega_u$ is a coset of the group
$$ E \coloneqq \{ F \in M(\mathrm{Y}, \T): dF \in \hat \Gamma \},$$
We introduce a countable dense sequence $G_1,G_2,\dots$ in $M(\mathrm{Y}, \T)$, and for each $u \in U$, let $n_u$ be the first integer such that there exists $F_u \in \Omega_u$ such that $\| e(F_u) - e(G_{n_u}) \|_{L^2(\mathrm{Y})} < \frac{1}{100}$; such an integer exists by density.  Assume for the moment that $n_u$ depends in a measurable fashion on $u$.  By \cite[Lemma C.1]{btz}, all the $F_u$ that arise in the above fashion differ from each other by a constant for fixed $u$.  In particular, there is a unique $F_u \in \Omega_u$ that minimizes $\| e(F_u) - e(G_{n_u}) \|_{L^2(\mathrm{Y})}$, and this $F_u$ clearly depends in a measurable fashion on $u$.  Setting $c_u \coloneqq h_u - dF_u$, we obtain Proposition \ref{mes-select} as claimed.

It remains to establish the measurability of $n_u$, which was asserted as being ``clearly'' true in \cite{btz}.  Clearly it suffices to show that for each $n$, the set
$$ \left\{ u \in U: \|e(F_u) - e(G_n)\|_{L^2(X)} < \frac{1}{100} \hbox{ for some } F_u \in \Omega_u \right\}$$
is measurable in $U$.  

Fix $n$.  Let $Z^1(\Gamma,\mathrm{Y},\T) \subset  M(\mathrm{Y},\T)^\Gamma$ denote the collection of $(\mathrm{Y},\T)$-cocycles.  The above set is the preimage under the map $u \mapsto h_u$ of the set
\begin{equation}\label{xs}
\{ h \in Z^1(\Gamma,\mathrm{Y},\T): h - dF \in \T^\Gamma \hbox{ for some } F \hbox{ with } e(F) \in B_n \}
\end{equation}
where
$$ B_n \coloneqq \left\{ f \in M(\mathrm{Y},S^1): \|f - e(G_n) \|_{L^2(X)} < \frac{1}{100} \right\}$$
(note that $h - dF$ is a cocycle, and so if it lies in $\T^\Gamma$ then it must come from a character in $\hat \Gamma$).
By the measurability of the map $u \mapsto h_u$, it suffices to show that \eqref{xs} is measurable in $Z^1(\Gamma,\mathrm{Y},\T)$.

The constraint $h - dF \in \T^\Gamma$ can be expanded as an equation of the form
$$ \frac{e(h_\gamma(x)) e(F(x))}{e(F(T^\gamma x))} = e(c_\gamma)$$
holding almost everywhere in $x$ for each $\gamma \in \Gamma$ and some $c_\gamma \in \T$.  If we now define the unitary operators $U^\gamma_h \colon L^2(\mathrm{Y}) \to L^2(\mathrm{Y})$ by
$$ U^\gamma_h f(x) \coloneqq e(-h_\gamma(x)) F(T^\gamma x)$$
(noting from the cocycle equation that these give a unitary action of $\Gamma$) and define a \emph{joint eigenfunction} of $(U_{\gamma,h})_{\gamma \in \Gamma}$ to be a function $f \in M(X,S^1)$ such that $U_{\gamma,h} f = \lambda_\gamma f$ holds for all $\gamma \in \Gamma$ and for some $\lambda_\gamma \in S^1$, we see that the set \eqref{xs} can be written as
$$ \{ h \in Z^1(\Gamma,\mathrm{Y},\T): (U_{\gamma,h})_{\gamma \in \Gamma} \hbox{ has a joint eigenfunction in } B_n\}.$$

For any $n'$ and sufficiently small $\eps>0$, we will show that there is a measurable set $S_{n',\eps}$ which contains
\begin{equation}\label{hmg}
 \{ h \in Z^1(\Gamma,\mathrm{Y},\T): (U_{\gamma,h})_{\gamma \in \Gamma} \hbox{ has a joint eigenfunction } f \hbox{ with } \|f-e(G_{n'}) \|_{L^2(\mathrm{Y})} < \eps^{10} \}
 \end{equation}
and is contained in
\begin{equation}\label{hmg-2}
 \{ h \in Z^1(\Gamma,\mathrm{Y},\T): (U_{\gamma,h})_{\gamma \in \Gamma} \hbox{ has a joint eigenfunction } f \hbox{ with } \|f-e(G_{n'}) \|_{L^2(\mathrm{Y})} < \eps \}
\end{equation}
taking a suitable countable union of such sets, we obtain the claim.

We now set $S_{n',\eps}$ to be the set
$$ S_{n',\eps} \coloneqq \{ h \in Z^1(\Gamma,\mathrm{Y},\T): \lim_{n \to \infty} \frac{1}{|\Phi_n|} \sum_{\gamma \in \Phi_n} |\langle U_{h}^\gamma G_{n'}, G_{n'} \rangle|^2 \geq 1-\eps^8 \},$$
where $\Phi_n$ is some F{\o}lner sequence for $\Gamma$.  The existence of the limit here follows from the mean ergodic theorem for Hilbert spaces (applied to the unitary action $\gamma \mapsto U_h^\gamma \otimes (U_h^\gamma)^*$ of $\Gamma$ on $L^2(\mathrm{Y}) \otimes L^2(\mathrm{Y})$).  Observe that $S_{n',\eps}$ is measurable.

Suppose that $h$ lies in the set \eqref{hmg}, then \footnote{Here we use the asymptotic notation $O(X)$ to denote a quantity bounded in magnitude by $CX$ for some absolute constant $C$.} 
$$\langle U_h^\gamma e(G_{n'}), e(G_{n'}) \rangle = \langle U_h^\gamma f, f \rangle + O(\eps^{10}) = 1 + O(\eps^{10})$$
for every $\gamma \in \Gamma$, and so $h \in S_{n',\eps}$ if $\eps$ is small enough.  Conversely, suppose that $h$ lies in the set $S_{n',\eps}$.  The operator
$$ Af \coloneqq \lim_{n \to \infty} \frac{1}{|\Phi_n|} \sum_{\gamma \in \Phi_n} \langle f, U_h^\gamma e(G_{n'}) \rangle U_h^\gamma e(G_{n'})$$
(with the limit existing in the weak operator topology at least, thanks to the mean ergodic theorem for Hilbert spaces as before) is a self-adjoint Hilbert--Schmidt operator of Hilbert--Schmidt norm at most $1$ (it is the limit finite rank operators of this form), and by construction one has
$$ \langle A e(G_{n'}), e(G_{n'}) \rangle \geq 1 - \eps^8.$$
From the spectral theorem, $A$ has a one-dimensional eigenspace of eigenvalue $1-O(\eps^8)$ (and all other eigenvalues of size at most $O(\eps^4)$, to maintain the Hilbert--Schmidt bound), and a unit eigenvector $f$ in this eigenspace is such that
$$ \langle e(G_{n'}), f \rangle \geq 1 - O(\eps^4)$$
and hence by the parallelogram law
$$ \| e(G_{n'}) - f \|_{L^2(\mathrm{Y})} = O(\eps^2).$$
Observe from the fact that the $U^\gamma_h$ are a group action and the F{\o}lner property that $A$ commutes with every $U_h^\gamma$, hence by one-dimensionality of the eigenspace, $f$ is a joint eigenfunction.  (Note that $|f|$ is $G$-invariant, hence constant by ergodicity, hence equal to $1$ since $f$ is a unit vector in $L^2(\mathrm{Y})$, so $f$ lies in $M(\mathrm{Y},S^1)$.)  Thus $h$ lies in the set \eqref{hmg-2}, and the claim follows.

\bibliographystyle{amsplain}

\end{document}